\documentclass[11pt]{amsart}
 \font\de=cmssi10
\usepackage{amssymb,amsmath}
\usepackage{amsthm}
\usepackage{amsfonts}
\usepackage[dvips]{graphicx}
\usepackage{psfrag}
\usepackage{bibentry}
\usepackage{enumerate}
\usepackage{euscript}
\usepackage{yfonts}

\newcommand{\F}{\mathcal{F}}
\newcommand{\T}{\mathbb{T}}

\renewcommand{\S}{\mathbb{S}}
\newcommand{\Pe}{\mathcal{P}}
\newcommand{\Res}{\mathcal{R}}
\newcommand{\A}{\mathcal{A}}
\newcommand{\B}{\mathcal{B}}

\newcommand{\R}{\mathbb{R}}

\newcommand{\Z}{\mathbb{Z}}

\newcommand{\He}{\mathcal H}
\newcommand{\W}{\mathcal W}
\newcommand{\x}{\mathbf x}
\newcommand{\s}{\sigma}
\newcommand{\go}[1]{\mathfrak {#1}}
\newcommand{\he}{\go{h}}
\newcommand{\g}{\go{g}}
\newcommand{\au}{\go{a}}

\renewcommand{\epsilon}{\varepsilon}

\newcommand{\vol}{\mathsf v}

\numberwithin{equation}{section}

\newtheorem{theorem}{Theorem}[section]

\newtheorem{corollary}[theorem]{Corollary}
\newtheorem{lema}[theorem]{Lemma}
\newtheorem{proposition}[theorem]{Proposition}

\newtheorem{conjecture}{Conjecture}

\newtheorem{problem}{\sc Problem}

\title{A survey on partially hyperbolic dynamics.}
\author{F. Rodriguez Hertz}
\author{M.A. Rodriguez Hertz}
\author{R. Ures}
\thanks{This work was partially supported by FCE 9021, CONICYT-PDT 29/220 and CONICYT-PDT 54/18 grants}
\address{IMERL-Facultad de Ingenier\'\i a\\ Universidad de la
Rep\'ublica\\ CC 30 Montevideo, Uruguay} \email{frhertz@fing.edu.uy}
\email{jana@fing.edu.uy} \email{ures@fing.edu.uy}
\date{\today}
\begin{document}

\begin{abstract}
Some of the guiding problems in partially hyperbolic systems are the
following: (1) Examples, (2) Properties of invariant foliations, (3)
Accessibility, (4) Ergodicity, (5) Lyapunov exponents, (6)
Integrability of central foliations, (7) Transitivity and (8)
Classification. Here we will survey the state of the art on these
subjects, and propose related problems.
\end{abstract}
\maketitle

{\footnotesize \tableofcontents}
\begin{section}{Introduction}
Here we will survey the state of the art in the area of partially
hyperbolic dynamics, that is, diffeomorphisms that display some
hyperbolic behavior ruling an intermediate one.\par
Hyperbolic behavior has proved to be a powerful tool to get
different types of chaotic properties from the ergodic and
topological viewpoints. As early as the late 60's or early 70's the
need of relaxing the full hyperbolicity hypothesis appeared. Indeed,
Pugh and Shub \cite{push1} in their study of the ergodicity of
Anosov actions and with Hirsch \cite{hipush1}, \cite{hipush} in
their study of invariant manifolds proposed the notion of normal
hyperbolicity, the intermediate part being played by a foliation
transversal to the hyperbolic part; Brin and Pesin \cite{brpe2},
\cite{brpe} studying ergodicity of skew products and frame flows,
proposed the notion of partial hyperbolicity where the intermediate
part is assumed to be tangent to a bundle. Both approaches are
obviously quite related and they essentially contain the same known
examples. Here we will follow this partially hyperbolic approach,
that is we will be assuming that $f:M\to M$ leaves a splitting
$TM=E^s\oplus E^c\oplus E^u$ invariant, where vectors in $E^s$ are
exponentially contracted in the future and vectors in $E^u$ are
exponentially contracted in the past. The dynamics in the
intermediate bundle $E^c$ has its point-wise spectrum between the
ones of $E^s$ and $E^u$.\par
In the last decade the area became quite active, here we will focus
on the following problems (1) Examples, (2) Properties of invariant
foliations, (3) Accessibility, (4) Ergodicity, (5) Lyapunov
exponents, (6) Integrability of central foliations, (7) Transitivity
and (8) Classification. We think that the theory will still be
growing, and that the formulation of problems, the basic and even
the simple ones, should be one of the main tasks.\par
Finally, we were not be able to cover all branches of study and some
of them were treated only laterally. Examples of these are Pesin
theory,  the search of SRB measures or Gibbs states,  the dominated
splitting approach proposed by Ma\~n\'e, partially hyperbolic
actions by more general groups, partially hyperbolic maps that are
not diffeomorphisms, etc. We encourage the reader also to read the
works \cite{bupushwi}, \cite{push4}, \cite{pe2}, \cite{bodivi} and
\cite{haspe}.

\par%
This work will be organized as follows. In section \ref{examples},
we enumerate the known examples of partially hyperbolic
diffeomorphisms. It is worth noting that this theory grows mainly
from examples. We propose different viewpoints to treat them. First,
in section \ref{section.invfol}, we present the invariant foliations
viewpoint. There we explore different kinds of regularities to be
found (or not) in the invariant foliations. In section
\ref{accessibility}, we focus on accessibility, that is, on the
ability to access from a point $x$ to another one $y$, by moving
only along lines piecewise tangent to $E^s$ and $E^u$, the bundles
composing the ruling hyperbolic part. Some results on accessibility
will be commented, and we shall compute accessibility classes on
three-dimensional nil-manifolds, and the general affine case.
Accessibility classes of stably ergodic toral automorphisms are
analyzed, and a short sketch of the proof that ergodic toral
automorphisms are stably ergodic is presented.\par There is evidence
that accessibility is abundant, and indeed there is a known
conjecture of Pugh and Shub that accessibility would be open and
dense; and also, that accessibility implies ergodicity among
partially hyperbolic diffeomorphisms. In sections
\ref{accessibility} and \ref{ergodicity}, we talk about advances
there have been in that sense, under certain restrictions on the
center bundle $E^c$. \par In section \ref{lyapunovexponent}, the
Lyapunov exponents viewpoint, and other kind of volume growth rates
is analyzed for partially hyperbolic systems.\par A completely
different scope is the study of the integrability of the center
distribution. This area seems quite open, and problems are posed.
This is reviewed in section \ref{integrability}. Section
\ref{rtransitivity} is devoted
 to robust transitivity, a point of view related to ergodicity.
 Finally, in section \ref{clasif}, the classification problem is
 reviewed.
\\
\emph{Acknowledgements.} The authors want to thank the Fields
Institute for warm hospitality and financial support during their
visit in January 2006.
\end{section}


\begin{section}{Examples}\label{examples}
Before going into the examples, let us define some terms more
precisely. A diffeomorphism $f:M\to M$ is {\de partially hyperbolic}
if there is an invariant splitting of the tangent bundle
$TM=E^s\oplus E^c\oplus E^u$ such that for all $x\in M$ and all
unitary vectors $v^\s\in E^\s$, with $\s=s,c,u$, the following
inequalities apply:
$$|D_xfv^s|<|D_xfv^c|<|D_xfv^u|$$
It is also required that $|Df|_{E^s}|<1$ and $|Df|_{E^u}|>1$.

As it will be mentioned below, there are two invariant foliations
$\F^s$ and $\F^u$, the {\de stable} and {\de unstable foliations},
that are tangent, respectively, to $E^s$ and $E^u$. These are the
only foliations with these features. But, in general, there is no
invariant foliation tangent to $E^c$; and, in case there were, it is
not known if it must be unique.\par
In studying partially hyperbolic systems, one of the problems is
that it is not clear if the amount of existing examples is small, or
if it essentially includes all the examples. Thus we get two
parallel problems: the search of examples and the classification
problem. We would like to split the examples into two categories in
nature, a grosser
 or topological one and another finer or geometric one; or even a
measure theoretic one.\par
For the topological type we would be interested in knowing in which
manifolds and in which homotopy classes the partially hyperbolic
dynamics can occur. For example, we say that two partially
hyperbolic systems $f:M\to M$ and $g:N\to N$, both having a central
foliation $\F$ are {\de centrally conjugated} or conjugated modulo
the central direction \cite{hipush} if there is a homeomorphism
$h:M\to N$ such that
\begin{enumerate}[i)]
\item $h\left(\F_f(x)\right)=\F_g\left(h(x)\right)$
\item $h\left(f\left(\F_f(x)\right)\right)=g\left(h\left(\F_f(x)\right)\right)$ or, which is equivalent, $\F_g\left(h(f(x))\right)=\F_g\left(g(h(x))\right)$
\end{enumerate}
It would be useful to classify partially hyperbolic systems modulo
central conjugacy. It would be interesting also to have an analogous
concept when the central distribution is not integrable.

Bellow we give a list of some of the existing examples. We hope that
we had put there most of them.

The second type of examples typically live within the first type and
will be appearing along the survey.
\begin{subsection}{Anosov Diffeomorphisms}\cite{an1}
A diffeomorphism $f:M\to M$ is an {\de Anosov diffeomorphism} if its
derivative $Df$ leaves the splitting $TM=E^s\oplus E^u$ invariant,
where $Df$ contracts vectors in $E^s$ exponentially fast, and $Df$
expands vectors in $E^u$ exponentially fast.\par
Anosov systems are the hallmark of hyperbolic and chaotic behaviors.
Nevertheless, they are far from being completely understood. For
example, the following problem is still open.
\begin{problem}\label{clasifanosov}\cite{sm}
Is every Anosov diffeomorphism conjugated to an infra-nil-manifold
automorphism?
\end{problem}
When the manifold underlying the dynamics is a nil-manifold or if
the unstable foliation has codimension one the answer is yes,
\cite{fr}, \cite{ma}, \cite{ne}. For expanding maps (when every
vector is expanded by the derivative) the answer also is yes, they
are always conjugated to infra-nil-manifold endomorphisms,
\cite{sh}, \cite{gro}.  It would be interesting to get analogous
results for partially hyperbolic diffeomorphisms, or at least to
have an answer to the following:
\begin{problem} {\cite{brbuiv}, Section \ref{clasif}.}
Let $f$ be a partially hyperbolic diffeomorphism on $\T^3$. Is it
true that its action in homology is partially hyperbolic?
\end{problem}
From the ergodic point of view, Anosov diffeomorphisms are very much
better understood.
\begin{theorem}\cite{an2}
Volume preserving Anosov diffeomorphisms are ergodic.
\end{theorem}

The partially hyperbolic systems share lots of their properties with
the Anosov systems. Let us describe some of those properties. There
are two invariant foliations $\F^s$ and $\F^u$ tangent to $E^s$ and
$E^u$. Both foliations have smooth leaves (as smooth as the
diffeomorphism), but the foliations themselves are not smooth a
priori. In fact, although there are some interesting cases where the
invariant foliations are smooth, the general case is that they are
rarely smooth \cite{an4}, \cite{haswi}. Thus, it became an
interesting problem to study the transversal regularity of these
foliations. For example, it turned out that the holonomies of these
foliations are {\de absolutely continuous} \cite{an2}, \cite{an3},
\cite{ansi}, \cite{si} i.e. we say that a map
$h:\Sigma_1\to\Sigma_2$ is absolutely continuous if it sends zero
measure sets into zero measure sets, see subsection
\ref{absolutecontinuity} for more details. The importance of
absolute continuity of the holonomies is that it implies that
Fubini's theorem is true for these foliations, that is, a measurable
set $A$ has zero measure if and only if for a.e. point $x$ in $M$
the intersection of $A$ with the leaf through $x$ has zero leaf-wise
measure. It is worth mentioning that in smooth ergodic theory, when
dealing with any kind of hyperbolicity, the smooth regularity of the
system is typically required to be at least
$C^{1+\mbox{\tiny{H\"older}}}$. In fact, in the $C^1$ category the
following is still unknown:
\begin{problem}
Are there examples of non ergodic volume preserving Anosov
diffeomorphisms?
\end{problem}

\end{subsection}

\begin{subsection}{Geodesic Flows}
Let $V$ be an $n$-dimensional manifold and let $g$ be a metric on
$V$. On $TV$ it is defined the {\de geodesic flow} as follows. Given
a point $x\in V$ and a vector $v\in T_xV$ there is a unique geodesic
$\gamma$ with $\gamma(0)=x$ and $\dot{\gamma}(0)=v$. For $t\in\R$ we
define $\phi_t(x,v)=\left(\gamma(t),\dot{\gamma}(t)\right)$. It
follows that $|\dot{\gamma}(t)|=|v|$ for every $t\in\R$ or more
precisely, $g_{\gamma(t)}\left(\dot{\gamma}(t)\right)=g_x(v)$ for
every $t\in\R$. Thus the geodesic flow preserves the vectors of a
given magnitude. Let $M=T_1V$ be the bundle of unit vectors tangent
to $V$ and let us restrict the geodesic flow to $M$. It turns out
that if the sectional curvature is negative then the geodesic flow
is in fact an Anosov flow \cite{an1}. Indeed, for every unit vector
$v$ in $M$, $T_vM$ may be identified with the orthogonal Jacobi
fields. Thus, if we call $E^s$ the set of orthogonal Jacobi fields
that are bounded for the future and $E^u$ the set of orthogonal
Jacobi fields that are bounded for the past, then negative sectional
curvature implies that $TM=E^s\oplus E^0\oplus E^u$ and the vectors
in $E^s$ are exponentially contracted in the future, $E^0$ is the
one dimensional space spanned by the vector-field defining the
geodesic flow and the vectors in $E^u$ are exponentially contracted
in the past.

The geodesic flow preserves a natural measure defined on $M$, the
Liouville measure $Liou$. Let us first define a one-form $\eta$ over
$TV$ as follows: if $\omega\in TV$ and $\chi\in T_\omega TV$ then we
define $\eta_{\omega}(\chi)$ as being $\omega\cdot d_\omega p(\chi)$
where $x=p(\omega)$, $p:TV\to V$ is the canonical projection. It
turns out that $d\eta$ is a symplectic 2-form on $TV$ and that the
geodesic flow preserves this symplectic form. Thus,
$L=d\eta\wedge\dots\wedge d\eta$ (n-times) is a $2n$-form. The
restriction of $L$ to $M$ is the $(2n-1)$-form defining $Liou$.\par
It was for the geodesic flows on surfaces of negative curvature that
Hopf \cite{ho} developed the machinery now called the Hopf argument
to prove ergodicity w.r.t. $Liou$ and the antecedent for Anosov
work. For general manifolds of negative sectional curvature Anosov
proved that the geodesic flow is ergodic w.r.t. $Liou$. In fact, he
proved more generally that $C^2$ volume preserving Anosov systems
are ergodic, thus, since being an Anosov flow is an open condition,
they form an open set of ergodic flows \cite{an1}, \cite{an2},
\cite{an3}, \cite{ansi}.

The time-one map of the geodesic flow on negative curvature, i.e.
$\phi_1$, is naturally a partially hyperbolic diffeomorphism. It was
not until 1992 that Grayson, Pugh and Shub, \cite{grpush} proved
that the time-one map of the geodesic flow on a surface of constant
negative curvature is a {\de stably ergodic} diffeomorphism, that
is, as in the Anosov case, their perturbations remain ergodic.
Later, Wilkinson proved the same result but for variable curvature
\cite{wi}.

There is also a the related topological question about robust
transitivity for partially hyperbolic systems that remains widely
open. In \cite{bodi} it is proven that close to the time-one map of
the geodesic flow on a negatively curved surface there are whole
open sets of transitive diffeomorphisms. But the following is still
open:
\begin{problem}
Is the time-one map of the geodesic flow on a negatively curved
surface robustly transitive?
\end{problem}
As it will be seen in section \ref{rtransitivity} it is enough to
prove that, for perturbations, the non-wandering set is still the
whole manifold.
\end{subsection}

\begin{subsection}{Anosov Flows}
We say that a flow $\phi_t$ on a manifold $M$ is an {\de Anosov
flow} if admits an invariant splitting $TM=E^s\oplus E^0\oplus E^u$,
where, as usual, vectors in $E^s$ and $E^u$ are respectively
contracted and expanded, and $E^0$ is the space spanned by the
vector-field. One of the main difference between Anosov flows and
Anosov diffeomorphisms is that there are known examples of Anosov
flows where the non-wandering set is not the whole manifold,
\cite{frwi}. In fact, this changes completely the hope of finding a
complete classification of Anosov flows like that stated in Problem
\ref{clasifanosov}. On the other hand, when dealing with transitive
Anosov flows, there is a dichotomy, either they are mixing or else
the bundle $E^s\oplus E^u$ is jointly integrable. In fact, in
\cite{pl}, it is proven that either the strong unstable manifold is
minimal or $E^s\oplus E^u$ is integrable. In this second case Plante
also proved that the flow is conjugated to a suspension but possibly
changing the time. In fact it is still an open problem to know if
the $su-$foliation is by compact leafs, and this is closely related
to the following long-standing problem
\begin{problem}
Is the action in homology of an Anosov diffeomorphism hyperbolic?
\end{problem}

As already mentioned, volume preserving Anosov flows are ergodic.
Moreover, the following is proven in \cite{bupuwi}:
\begin{theorem}
Let $\phi_1$ be the time-one map of a volume preserving Anosov flow
$\phi$. Then $\phi_1$ is stably ergodic if and only if $\phi$ is
mixing.
\end{theorem}
Of course, the time-one map of the suspension of an Anosov
diffeomorphism by a constant roof function is not stably ergodic.
\end{subsection}

\begin{subsection}{Frame Flows}\cite{br1}, \cite{brgr}, \cite{brka}, \cite{brpe2}, \cite{bupo}.
The {\de frame flow} on a riemannian manifold $(V,g)$ fibers over
its geodesic flow. Let $\hat{M}$ be the space of positively oriented
orthonormal $n$-frames in $TV$. Thus $\hat{M}$ naturally fibers over
$M=T_1V$, where the projection takes a frame to its first vector.
The associated structure group $SO(n-1)$ acts on fibres by rotating
the frames keeping the first vector fixed. In particular, we can
identify each fiber with $SO(n-1)$. Let
$\hat{\phi}_t:\hat{M}\to\hat{M}$ denote the frame flow, which acts
on frames by moving their first vectors according to the geodesic
flow and moving the other vectors by parallel transport along the
geodesic defined by the first vector. The projection is a
semi-conjugacy from $\hat{\phi}_t$ to $\phi_t$. In particular,
$\hat{\phi}_t$ is an $SO(n-1)$-group extension of $\phi_t$. The
frame flow preserves the measure $\mu=Liou\times\nu_{SO(n-1)}$,
where $\nu_{SO(n-1)}$ is the (normalized) Haar measure on $SO(n-1)$.
It turns out that the time-$t$ map of the frame flow is a partially
hyperbolic diffeomorphism \cite{brpe}. The neutral direction has
dimension $1+dimSO(n-1)$ and is spanned by the flow direction and
the fibre direction.\par
The frame flow on manifolds of negative sectional curvature is known
to be ergodic in lots of cases. The study of the ergodicity of the
frame flow restricts to the study of its accessibility classes (see
section \ref{accessibility} for the notion of accessibility) and is
a very interesting example to begin with, in order to learn how to
manage them. Finally the frame flow is stably ergodic in the cases
it is known to be ergodic. But it is not always ergodic, K\"ahler
manifolds with negative curvature and real dimension at least 4 have
non-ergodic frame flows because the complex structure is invariant
under parallel translation. We suggest the reader to see \cite{bupo}
for a good account of the existing results, problems  and
conjectures.

\end{subsection}



\begin{subsection}{Affine diffeomorphisms}\label{afdif}
Let $G$ be a Lie group and $B\subset G$ a subgroup. Given a one
parameter subgroup of $G$ it defines an {\de homogeneous flow} on
$G/B$. Examples of homogeneous flows are  geodesic flows of
hyperbolic surfaces. There are lots of interplays between the
dynamics of homogeneous flows and the algebraic properties of the
groups involving it, see for example \cite{sta} for an account.

The time-$t$ map of an homogeneous flow is a particular case of an
{\de affine diffeomorphism}. In fact affine diffeomorphisms and
homogeneous flows are typically treated in a similar way. Let $G$ be
a connected Lie group, $A:G\to G$ an automorphism, $B$ a closed
subgroup of $G$ with $A(B) = B$, and $g\in G$. Then we define the
affine diffeomorphism $f:G/B\to G/B$ as $f(xB)=gA(x)B$. We shall
assume that $G/B$ supports a finite left $G-$invariant measure and
call, in this case, $G/B$ a finite volume homogeneous space. If
$G/B$ is compact and $B$ is discrete the existence of such a measure
is immediate, but if $B$ is not discrete the assumption is
nontrivial.

The affine diffeomorphism $f$ is covered by the diffeomorphism
$\bar{f}=L_g\circ A:G\to G$; where $L_g:G\to G$ the left
multiplication by $g$. If $\g$ is the Lie algebra of $G$, we may
identify $T_eG=\g$ where $e$ is the identity map. Let us fix a right
invariant metric on $G$, i.e. $R_g$ is an isometry for every $g$
where $R_g$ is right multiplication by $g$. Let us define the
naturally associated automorphism $\au(f):\g\to\g$ by
$\au(f)=Ad(g)\circ D_eA$ where $Ad(g)$ is the adjoint automorphism
of $g$, that is the derivative at $e$ of $x\to gxg^{-1}$. In other
words, $\au(f)$ is essentially the derivative of $\bar{f}$, but
after right multiplication by $g^{-1}$ (which is an isometry) in
order to send $T_gG$ to $T_eG$. So we have the splitting
$\g=\g^s\oplus\g^c\oplus\g^u$ w.r.t the eigenvalues of $\au(f)$
being of modulus less than one, one, or bigger than one respectively
and similarly, $\g^s$ is formed by the vectors going exponentially
to $0$ in the future, $\g^u$ is formed by the vectors going
exponentially to $0$ in the past and $\g^c$ is formed by the vectors
that grow at most polynomially for the future and the past. Observe
that if $v_{\lambda}$ and $v_{\sigma}$ are eigenvectors for $\au(f)$
w.r.t. $\lambda$ and $\sigma$ respectively then we have that
$$
\au(f)\left([v_{\lambda},v_{\sigma}]\right)=
\left[\au(f)(v_{\lambda}),\au(f)(v_{\sigma})\right]=\lambda\sigma[v_{\lambda},v_{\sigma}]
$$
and hence if $[v_{\lambda},v_{\sigma}]\neq 0$ then it is an
eigenvector for $\lambda\sigma$. As a consequence we get that
$\g^s$, $\g^u$, $\g^c$, $\g^{cs}=\g^c\oplus\g^s$ and
$\g^{cu}=\g^c\oplus\g^u$ are subalgebras tangent to connected
subgroups $G^s, G^u, G^c, G^{cs}$ and $G^{cu}$ of $G$ and their
translates will define the stable, unstable, center, center-stable
and center-unstable foliations respectively.

Let $\he$ denote the smallest Lie subalgebra of $\g$ containing
$\g^s$ and $\g^u$. Using Jacobi identity it is not hard to see that
it is an ideal, $\he$, called the {\de hyperbolic subalgebra} of
$\bar{f}$. Moreover, let us denote $H\subset G$ the connected
subgroup tangent to $\he$ and call it the hyperbolic subgroup of
$\bar{f}$. As $\he$ is an ideal in $\g$, $H$ is a normal subgroup of
$G$. Finally let us denote with $\go{b}\subset\g$ the Lie algebra of
$B\subset G$. Then we have the following:
\begin{theorem}\cite{push4}
Let $f:G/B\to G/B$ be an affine diffeomorphism as above, then $f$ is
partially hyperbolic if and only if $\he\not\subset\go{b}$.
Moreover, if $f$ is partially hyperbolic then the left action of
$G^{\sigma}$, $\sigma=s, u, c, cs, cu$ on $G/B$ foliates $G/B$ into
the stable, unstable, center, center-stable and center-unstable
foliations respectively.
\end{theorem}

\begin{problem}
Is there an example of a non-Anosov affine diffeomorphism that is
robustly transitive? Are they exactly the same as the stably ergodic
ones?.
\end{problem}

\end{subsection}

\begin{subsection}{Linear Automorphisms on Tori}
A special case of affine diffeomorphisms are the affine
automorphisms on tori. In fact, the torus $\T^N$ may be seen as the
quotient $\R^N/\Z^N$. Integer entry $N\times N$ matrices with
determinant $\pm 1$ define what we shall call linear automorphisms
of tori simply via matrix multiplication. Thus, given such a matrix
$A$ and a vector $v\in\R^N$, it is defined an affine diffeomorphism
of the torus $f$ by $f(x)=Ax+v$. It is quite easy to see that,
conjugating by a translation, it is enough to study the case where
$v$ belongs to the eigenspace corresponding to the eigenvalue $1$,
$E_1$. Observe also that $E_1$ is a rational space, that is, it has
a basis formed by vectors of rational coordinates.\par
The corresponding splitting of the tangent bundle here, is the
splitting given by the eigenspaces of $A$. Thus, a not quite
involved argument proves that $f$ is partially hyperbolic unless all
the eigenvalues of $A$ are roots of unity. Moreover, using a little
bit of harmonic analysis (Fourier series) it is seen \cite{hal} that
$f$ is ergodic if and only if $A$ has no eigenvalues that are roots
of the identity other than one itself and $v$ has irrational slope
inside $E_1$. Finally, notice that if $E_1$ is not trivial, we may
always perturb in order to make $v$ of rational slope, thus in order
to get that perturbations remain ergodic it is necessary that also
$1$ be not in the spectrum of $A$. Thus we reach to the following
problem:
\begin{problem}\label{ergodicautomorphismproblem}\cite{hipush}, \cite{frh}, Subsection
\ref{ergodicautomorphism}. Are the ergodic linear automorphisms
stably ergodic?
\end{problem}
Of course, an analogous problem may be posed in the topological
category, that is, are their perturbations also transitive?
\cite{hipush}.
\end{subsection}
\begin{subsection}{Direct Products}
Given a partially hyperbolic diffeomorphism $f:M\to M$ and $g:N\to
N$ a diffeomorphism, the product $f\times g:M\times N\to M\times N$
is partially hyperbolic if the dynamics of $g$ is less expanding and
contracting, respectively, than the expansions and contractions of
$f$. This is essentially the most trivial way a partially hyperbolic
dynamics appears, Anosov$\times$identity. Besides, we can also make
the product of two partially hyperbolic diffeomorphisms.

It is quite interesting that by making perturbations of this product
dynamics, lots of nontrivial examples arises. For instance, the
first example of a robustly transitive non-Anosov diffeomorphism
constructed by Shub \cite{sh2}, see section \ref{rtransitivity},
although not a product, is a large perturbation of a product. In
fact direct products as well as the construction of Shub are part of
a more general type of construction, the partially hyperbolic
systems that fiber over other partially hyperbolic systems.
\end{subsection}

\begin{subsection}{Fiberings over partially hyperbolic diffeomorphisms.}
Let $f:B\to B$ be a partially hyperbolic diffeomorphism with
splitting $TM=E^s_f\oplus E^c_f\oplus E^u_f$. Let $p:N\to B$ be a
fibration with fiber $F$, let us call $F(x)$ the fiber through $x$.
Then any lift $g:N\to N$ of $f$ is a partially hyperbolic
diffeomorphism if
$$
\left|D_{p(x)}f|E^s_f\right|<m\left(D_xg|T_xF(x)\right)\leq\left|D_xg|T_xF(x)\right|<m\left(D_{p(x)}f|E^s_f\right).
$$
where $m(A)=|A^{-1}|^{-1}$. As we said, Shub's example of a robustly
transitive diffeomorphism is of this kind, and, in fact, many of the
existing examples are of this kind. It would be interesting to find
the minimal pieces over which partially hyperbolic systems are
built. For example:
\begin{problem}
Find the partially hyperbolic diffeomorphisms $f$ such that no
partially hyperbolic diffeomorphism $g$ homotopic to $f^n$, $n>0$,
fibers over a lower dimensional partially hyperbolic diffeomorphism.
The geodesic flow on negative curvature as well as the ergodic
automorphisms of tori defined in \cite{frh} are examples of that
building blocks. Find other types of gluing technics to generate new
partially hyperbolic systems.
\end{problem}
\end{subsection}

\begin{subsection}{Skew products.}
Another type of systems that fiber over lower dimensional partially
hyperbolic diffeomorphisms are the {\de skew products}: Let $f:M\to
M$ be a partially hyperbolic diffeomorphism, $G$ a Lie group and
$\theta:M\to G$ a function. Define the skew product
$f_{\theta}:M\times G\to M\times G$ by
$f_{\theta}(x,g)=(f(x),\theta(x)g)$. Skew products where extensively
studied in the context of partially hyperbolic diffeomorphisms, see
for example \cite{adkish}, \cite{br1}, \cite{br2}, \cite{brpe},
\cite{buwi1}, \cite{fipa}.

\end{subsection}

\end{section}

\begin{section}{Properties of the invariant
foliations}\label{section.invfol}
 By a foliation we shall mean a
topological foliation. In the dynamical framework it is useful to
treat the regularity of the leaves of a foliation and the regularity
of their holonomies separately. In fact, typically the foliations we
deal with have smooth leaves but the holonomies are only H\"older
continuous.\par
Given a continuous plane field $E\subset TM$, we say that $E$ is
{\de integrable} if there is a foliation $\F$ tangent to $E$, i.e.
each leaf of the foliation is a $C^1$ manifold everywhere tangent to
$E$.\par
Given a partially hyperbolic diffeomorphism $f:M\to M$, with
splitting $TM=E^s\oplus E^c\oplus E^u$, we are mainly interested in
studying if the plane fields are integrable, and in case they are,
how good their leaf and holonomy regularity are.

The study of the regularity of the invariant bundles and the
invariant foliations is in the core of the theory. There is a vast
bibliography on the subject and we will not be able to mention every
reference here. Nevertheless, the book of Hirsch-Pugh-Shub
\cite{hipush} is still one of the more inspiring ones in the
literature. So that this subject is quite well developed although it
is still very active and there continue to appear new types of
regularity to deal with. In fact one of the themes here is what we
will mean by regularity.

\begin{subsection}{Existence and regularity of the strong foliations}
First we shall deal with the case of the strong bundles, that is,
the stable and unstable bundles. Let us first go into the
integrability problem:
\begin{theorem}\cite{brpe}, \cite{hipush}
Given a partially hyperbolic diffeomorphism the stable and unstable
bundles are always integrable to foliations $\F^s$ and $\F^u$ whose
leaves are as differentiable as the diffeomorphism.
\end{theorem}
From this, we have that the integrability of the strong bundles is
quite well solved and understood. So, let us go into the problem of
the regularity of the holonomies. Let us split this problem. We
shall treat the regularity of the holonomies from a purely
differentiable viewpoint. In the next subsection we shall deal with
the absolute continuity problem. Given $r\in\Z$, $r\geq 0$ and
$0<\theta\leq 1$, we say that a map between smooth manifolds {\de is
$C^{r+\theta}$} if it is $C^r$ and its $r^{th}$ derivative is
H\"older continuous with H\"older exponent $\theta$. Of course, in
order to say that a map is $C^{r+\theta}$ we need first that the
manifolds where it is defined are at least $C^{r+\theta}$. This last
observation, though it is an obvious remark, becomes a real problem
when trying to prove regularity of the holonomies.\par
The regularity of the holonomies is much related to what one could
call the spectrum of the differential of  the partially hyperbolic
diffeomorphism. Instead of going into the definition of the Mather
spectrum (the reader may found a good exposition in \cite{haspe}) we
shall go into its point-wise version.

For a linear transformation $A$ between Banach spaces, we define
$m(A)=\inf_{|v|=1}|Av|$. Observe that if $A$ is invertible then
$m(A)=|A^{-1}|^{-1}$. \par
The following theorem is essentially proved in \cite{pushwi1}, we
thank Keith Burns for pointing out this statement.
\begin{theorem}\label{bunchimpliesdiff}
Let $M$ be a complete riemannian manifold and let $f:M\to M$ be a
$C^k$ diffeomorphism with an invariant splitting $TM=E_1\oplus E_2$
satisfying
$$
\sup_x\frac{|D_xf|_{E_1}|}{m(D_xf|_{E_2})}<1.
$$
Let us assume that
$$
\sup_x|D_xf|_{E_1}|\frac{|D_xf|_{E_2}|^r}{m(D_xf|_{E_2})}<1
$$
Then there is a foliation tangent to $E_1$ whose leaves are as
smooth as $f$ and whose holonomies are $C^l$ where
$l=\min\{k-1,r\}$. Here $r$ and $k$ are allowed to be any positive
real numbers bigger than or equal to $1$.
\end{theorem}
For example, if $f$ is a partially hyperbolic with $E^{cs}$
integrable and
\begin{eqnarray}\label{bunchparadiff}
\sup_x|D_xf|_{E^s}|\frac{|D_xf|_{E^c}|^r}{m(D_xf|_{E^c})}<1
\end{eqnarray}
then we have that the stable foliation is $C^l$ smooth when
restricted to each center-stable leaf, where $l$ depends on the
differentiability of $f$, the differentiability of the leaves of the
center-stable foliation, and $r$. But typically the leaves of the
center-stable foliation are not better than
$C^{1+\mbox{\tiny{H\"older}}}$. Nevertheless, if $f$ is $C^2$, the
stable foliation is still $C^1$ when restricted to each
center-stable leaf if equation \ref{bunchparadiff} holds with $r=1$,
\cite{pushwi1}. Moreover, it can be seen that if the center
dimension is one, then the stable foliation is still $C^1$ when
restricted to any $c+s$ dimensional manifold everywhere tangent to
$E^{cs}$ \cite{rhrhur1}. Notice that when the center dimension is
one, equation \ref{bunchparadiff} with $r=1$ is trivially satisfied.
\begin{problem}
Prove that the stable foliation is still $C^1$ when restricted to
any $(c+s)$-dimensional manifold everywhere tangent to $E^{cs}$
whenever equation (\ref{bunchparadiff}) is satisfied with $r=1$.
\end{problem}

Another interesting case is in dimension $3$, if $f$ is volume
preserving, then equation (\ref{bunchparadiff}) is satisfied with
$r=2$. But the problem is that a priori it is not known if there is
a $C^2$ center-stable manifold and hence theorem
\ref{bunchimpliesdiff} cannot be straightforwardly applied. It would
be interesting to know if the results in \cite{pushwi1},
\cite{pushwi2} can be adapted to solve the following:
\begin{problem}\label{difofsuhol}
If $f:M\to M$ is a volume preserving, partially hyperbolic $C^r$
diffeomorphism, $r$ big enough, with $dim M=3$, prove that the
stable and unstable holonomies are $C^2$ when restricted to some
suitable $2$ dimensional manifolds. At least prove that if
$E^s\oplus E^u$ is integrable then it integrates to a $C^2$
foliation.
\end{problem}

Related to this we want to ask the following:
\begin{problem}
What is the relation between bunching and the smoothness of
$E^s\oplus E^u$ if any?
\end{problem}
There are results when $f$ is only $C^{1+\mbox{\tiny{H\"older}}}$,
see \cite{buwi4}. Also, when there is no integrability of the
center-stable distribution at all, we still have the
differentiability of the stable foliation restricted to some ``fake"
foliations, if equation (\ref{bunchparadiff}) holds with $r=1$, see
for instance \cite{buwi3}. The reader may found problems related to
this in subsection \ref{smoothcase}.\par
Another related problem about the holonomies is how they vary as we
perturb the dynamics. In \cite{frh} it is proven that their
variation is quite good in the particular case analyzed there. It
would be interesting to have more general results with explicit
constants.

There is another form of regularity that is between the smooth
category and the absolute continuity category. It is the {\de
quasi-conformality} and it proved to be a very powerful notion at
the time of establishing ergodicity. We shall return to this subject
in subsection \ref{juliennes}.
\end{subsection}

\begin{subsection}{Absolute continuity of the strong invariant foliations}\label{absolutecontinuity}
Let us enter now into absolute continuity. To begin with, there are
several notions of absolute continuity of foliations. It seems that
the strong foliations satisfy the strongest one. But let us begin
with the weakest definition.\par
Given $\F$ a foliation let us denote with $W(x)$ the leaf of $\F$
through $x$. Given a measure $\mu$ we call $\mu_{W(x)}$ the
conditional measures along $W(x)$. We say that $\F$ is {\de
absolutely continuous with respect to $\mu$} if given a measurable
set $A$ we have that $\mu(A)=0$ if and only if there is a set $B$
such that $\mu(B)=1$ and such that $\mu_{W(x)}(A\cap W(x))=0$ for
every $x\in B$. Typically, the absolute continuity problem is
analyzed w.r.t. Lebesgue measure, but it would be interesting to
know the absolute continuity of the strong foliations with respect
to other measures. For example, it seems likely that the unstable
foliation is absolutely continuous with respect to any $u$-Gibbs
state.

We say that the foliation $\F$ {\de has absolutely continuous
holonomies} if given $x$, $y\in W(x)$ and two transversals
$\Sigma_x$ and $\Sigma_y$ then the holonomy map
$h:\Sigma_x\to\Sigma_y$ sends zero sets into zero sets. Again here
typically we analyze the case when the measure on the sections are
Lebesgue, but one may ask about other transversal measures.
Moreover, if $\F$ has absolutely continuous holonomies, then given
an holonomy $h:\Sigma_x\to\Sigma_y$ one may talk about the
Radon-Nikodym derivative of the pull back measure on $\Sigma_y$ over
the measure on $\Sigma_x$, and call it the {\de Jacobian} of $h$,
that is, if we denote with $\mu_{\Sigma_x}$ and $\mu_{\Sigma_y}$ the
measures on the transversals, then the jacobian will be the map
$Jh:\Sigma_x\to\R$ satisfying
$$
\mu_{\Sigma_y}(A)=\int_AJh(t)d\mu_{\Sigma_x}(t)
$$
for every measurable set $A\subset\Sigma_x$. So the absolute
continuity theorem is:
\begin{theorem}\label{abscontstrongfol}\cite{an1}, \cite{brpe},
\cite{push1} The stable and unstable foliation for a partially
hyperbolic diffeomorphism have absolutely continuous holonomies.
Moreover, the holonomies are uniformly H\"older continuous, that is,
the exponent and constant may be taken uniformly on the whole
manifold if the transversals are taken good enough (uniformly
smooth, bounded and reasonable angle with the foliation) and close
enough (distance between the points not going to infinity along the
leaf).
\end{theorem}

\end{subsection}

\begin{subsection}{Non-absolutely continuity of central foliations}\label{fubinineightmare}

In this section we shall see how absolute continuity fails
completely for the central foliation. In fact, there are examples
with a full measure set that intersects each leaf only in a finite
number of points. This phenomenon is sometimes known as Fubini's
nightmare.

In \cite{mi} the reader may found the first example, built by
Anatole Katok, of such phenomenon of a non-absolutely continuous
foliation. The idea is essentially as follows, take a smooth path
$f_t$, $t\in [0,1]$ of area preserving Anosov diffeomorphisms on
$\T^2$ beginning with a linear one. Then one can define a partially
hyperbolic diffeomorphism $F:\T^2\times [0,1]\to \T^2\times
[0,1]\to$ by $F(x,t)=(f_t(x),t)$. It follows that $F$ leaves a
center foliation invariant, and that the leaves through $(x,0)$ can
be parameterized by $(h_t(x),t)$ where $h_t:\T^2\to\T^2$, $t\in
[0,1]$ is a conjugacy homotopic to identity between $f_0$ and $f_t$,
i.e. $f_t\circ h_t=h_t\circ f$. Hence we have that $h_t$ is the
central holonomy between the transversal $\T^2\times\{0\}$ and
$\T^2\times\{t\}$. Now, if we define
$\mu_t(A)=Leb\left(h_t^{-1}(A)\right)$, then $\mu_t$ is the entropy
maximizing measure for $f_t$. On the other hand it is known that,
for an Anosov diffeomorphism of $\T^2$, if $Leb$ is the entropy
maximizing measure then the eigenvalues at periodic points for the
diffeomorphism coincide with the ones of the linear model, moreover,
it is smoothly conjugated to the linear model. So, if we take the
path $f_t$ in such a way that the eigenvalues at the fixed point are
different from the linear one for every $t$, then $\mu_t$ is not
Lebesgue measure and hence it must be singular. Hence, for every
$t>0$, there is a set $A_t$ of full Lebesgue measure on $\T^2$ whose
image under the center holonomy $h_t(A_t)$ has zero Lebesgue
measure. In fact, if the eigenvalue at the fixed point of any two
$f_t$ are different, then there is a full measure set
$B\subset\T^2\times [0,1]$ that cuts each central leaf at exactly
one point. Take $B$ to be the set of points $(x,t)$ such that its
time average converges to its space average w.r.t. $Leb_t$, the
Lebesgue measure on $\T^2\times\{t\}$. Then $B$ has full measure
because it has full measure when restricted to each torus
$\T^2\times\{t\}$. Now, given $t\neq s$, we can send $Leb_t$ to
$\T^2\times\{s\}$ by central holonomy (that correspond to
conjugating $f_t$ with $f_s$) and get a measure that should be
singular w.r.t. $Leb_s$ (this follows from a theorem of de la LLave
\cite{dL}). Hence, if $(x,t)\in B$,  the corresponding point $(y,s)$
by central holonomy would be typical for the measure singular to
$Leb_s$ and hence $(y,s)\notin B$, in other words, $B$ intersects
each central leaf at exactly one point.\par
In \cite{shwi2}, Mike Shub and Amie Wilkinson found the same
phenomenon with a different approach, thus finding an open set where
the center foliation is not absolutely continuous. In this case, $F$
is a diffeomorphism of $\T^3=\T^2\times\S^1$ close to a skew product
over an Anosov diffeomorphism of $\T^2$. They built an ergodic
volume preserving diffeomorphism with nonzero central exponents, see
section \ref{lyapunovexponent}. Since all central curves are compact
it follows from an argument by Ma\~n\'e that the central foliation
cannot be absolutely continuous. Indeed the argument is as follows:
assume that the central exponent is positive. If the foliation were
absolutely continuous then, using Pesin theory, one can find a
positive leaf measure set where there is actual central expansion,
i.e. there is a set $A\subset W^c$ of positive length measure such
that $|D_xF^n|E^c|>\sigma^n$, $\sigma>1$ for $x\in A$, and every
$n\geq n_0$. But then
$$
C\geq length\left(f^n(W^c)\right)\geq length\left(f^n(A)\right)\geq
\sigma^n length(A)
$$
a contradiction. Here it is also possible to find a full measure set
$B\subset\T^3$ that cuts each central curve in a finite set, see
\cite{ruwi}.

In \cite{hipe}, Hirayama and Pesin generalize Ma\~n\'e's argument
and prove the following general theorem.

\begin{theorem}
Let $f$ be a $C^2$ diffeomorphism of a compact smooth Riemannian
manifold $M$ preserving a smooth measure $\mu$. Let also $\W$ be an
$f$-invariant foliation of $M$ with smooth leaves. Assume that $\W$
has finite volume leaves almost everywhere. If $f$ is
$\W$-dissipative then the foliation $\W$ is not absolutely
continuous.
\end{theorem}
Here, $f$ is said to be {\de $\W$-dissipative} if
$$
\int_M \log Jac\left(Df|\W\right)d\mu\neq 0.
$$

Another type of example is as follows, take a linear Anosov
diffeomorphism $A$ of $\T^3$ with three one-dimensional invariant
bundles, $E^{uu}\oplus E^u\oplus E^{ss}$. Let us call $\chi^{uu}$,
$\chi^u$ and $\chi^{ss}$ the corresponding Lyapunov exponents of $A$
and $\chi^{uu}_f$, $\chi^u_f$ and $\chi^{ss}_f$ the ones w.r.t.
Lebesgue measure of a volume preserving perturbation $f$. The first
observation is that if $f$ is a volume preserving perturbation of
$A$, since the strong foliations are absolutely continuous, it
follows that the strong stable exponent cannot decrease, and the
strong unstable Lyapunov exponent cannot increase, i.e.
$-\chi^{ss}_f\leq -\chi^{ss}$ and $\chi^{uu}_f\leq\chi^{uu}$. On the
other hand, if we make the perturbation in such a way that it
preserves the $E^{uu}\oplus E^u$ foliation then the strong stable
exponent will be preserved and hence the entropy (w.r.t. Lebesgue
measure) will not change, notice that the entropy also equals
$\chi^{uu}_f+\chi^u_f$. Similarly to the strong case, if the central
foliation were absolutely continuous then the central exponent will
not increase either, i.e. $\chi^u_f\leq \chi^u$. Then the final
stage to get the non absolutely continuous central foliation will be
to make the perturbation so that $\chi^{uu}_f<\chi^{uu}$. One way to
make this perturbation is along the lines of \cite{shwi2}. Another
way is just to notice that if $\chi^{uu}=\chi^{uu}_f$, then the
strong unstable eigenvalues at any periodic point coincides with the
linear ones, so brake any of them and we are done.


We found this example in a conversation with Anatole Katok, Andrey
Gogoleg, a student of Katok, and the first author. One week later,
in the International Workshop on Global Dynamics Beyond Uniform
Hyperbolicity at Chicago, in his talk, Radu Saghin exposed his work
with Zhihong Xia about the construction of similar examples within
the study of different kinds of dynamical growths (see section
\ref{lyapunovexponent}).

It is generally believed that the failure of absolute continuity of
the central foliation is a generic phenomenon.
\begin{problem}
Prove that there is a $C^1$ open, $C^{\infty}$ dense, set of
diffeomorphism whose central foliation fails to be absolutely
continuous. Moreover, try to characterize the case where the central
foliation is absolutely continuous, at least when the central
dimension is one.
\end{problem}

We finish with another problem.
\begin{problem}
Analyze regularity of the central foliation, maybe when restricted
to the center-stable or center-unstable. For example, if the
dynamics on the central direction is an isometry, is there some type
of regularity of the the central foliation?. What about the case
when there are no periodic points or when all the central Lyapunov
exponents are $0$.
\end{problem}
Observe that in Katok's original example, the dynamics on the
central direction is an isometry for some suitable metric. In fact,
if we take the sup norm $|v|=\max\{|v_{\T_2}|,|v_{\S^1}|\}$ then it
follows that $|D_xF|E^c|=1$. Hence absolute continuity of the
central foliation seems to imply some kind of strong rigidity
phenomenon, not only on the growth of the central direction but also
on the strong directions. It would be interesting to find out what
can be said when the central foliation is absolutely continuous. For
example, is there some kind of converse to Katok's original example
when the diffeomorphism is a perturbation of Anosov times identity?
and finally,
\begin{problem}
Can we have zero-central exponent, accessibility and non-absolutely
continuous central foliation?, what about for perturbations of
Anosov times identity on $\T^2\times\S^1$?
\end{problem}
\end{subsection}

\begin{subsection}{}
\end{subsection}

\end{section}

\begin{section}{Accessibility}\label{accessibility}
Given two uniquely integrable sub-bundles $E,F\subset TM$, an
equivalence relation may be defined by saying that $y$ is {\de
$(E,F)-$accessible} from $x$ if there is a piecewise smooth path,
piecewise tangent to either $E$ or $F$ beginning at $x$ and ending
at $y$. Another way to define it is to say that the accessibility
class of $x$ is the minimal set containing $x$, saturated by leaves
of the $E$-foliation and the $F$-foliation. When $(E,F)=(E^s,E^u)$
corresponding to a partially hyperbolic system, we call it the {\de
$su$-accessibility relation} and we denote the $su-$accessibility
class of $x$ by $C(x)$. The study of the accessibility classes is
quite well developed in control theory, but typically, in control
theory the bundles are assume to be smooth. In dynamics, as we
already have said, the bundles are rarely smooth, so that much of
the work already done in control theory should be redone in this new
setting. Some properties follow straightforward, but others become
quite difficult.\par
To the best of our knowledge, the use of the accessibility property
to prove ergodicity was first used by Sacksteder in \cite{sa}. He
essentially proved that accessibility implies ergodicity when the
strong foliations are smooth. Later, Brin and Pesin in \cite{brpe}
used it again in the context of frame flows and skew products where
they also look at its relation with transitivity, see subsection
\ref{weakergodic}.
Finally, Pugh and Shub, beginning with
\cite{grpush}, used it systematically within a plan to prove
ergodicity for partially hyperbolic systems. Indeed, Pugh and Shub
asked the following:
\begin{conjecture} \label{mainconjectureaboutergodicity}
Stable ergodicity is $C^r$-dense among partially hyperbolic system,
$r\geq 2$.
\end{conjecture}
Their plan is to split it into two conjectures related to the
accessibility.
\begin{conjecture}\label{accessimplieserg}
Essential accessibility implies ergodicity
\end{conjecture}
Here {\de essential accessibility} means that any measurable
$su$-saturated set has either full or null measure.
\begin{conjecture}\label{accesisgeneric}
Stable accessibility is $C^r$-dense among partially hyperbolic
systems, volume preserving or not, $r\geq 2$.
\end{conjecture}
When $\dim E^c=1$ conjecture \ref{accesisgeneric} is proven in its
full strength in \cite{rhrhur1} for the volume preserving case and
in \cite{burhrhtaur} for the non preserving   case, see also
subsection \ref{accescentral1d} for an account on this. Moreover, in
\cite{rhrhur1}, the authors also prove conjecture
\ref{accessimplieserg} when $\dim E^c=1$, thus proving the main
conjecture \ref{mainconjectureaboutergodicity} when $\dim E^c=1$. We
removed the dynamical coherence hypothesis of the paper
\cite{buwi2}, when $\dim E^c=1$, essentially using the notion of
weak-integrability in \cite{brbuiv}. At the same time, Keith Burns
and Amie Wilkinson also removed the dynamical coherence hypothesis
for any central dimension (under a bunching assumption that is
trivially satisfied if $\dim E^c=1$) but using the idea of fake
foliations, \cite{buwi3}, see subsection \ref{ergodicity} where
conjecture \ref{accessimplieserg} will be treated.\par
Conjecture \ref{accesisgeneric} is also known to be true in full
generality but only for $r=1$ by \cite{dowi}. There are also lot of
special cases where conjecture \ref{accesisgeneric} holds, that is,
systems that are known to be stably accessible or that can be
approached by stably accessible ones. Historically, the first non
Anosov examples having the accessibility property were the ones in
Sacksteder work, \cite{sa}, for some affine diffeomorphisms. Then
Brin \cite{br1}, \cite{br2} in his work on skew products proved that
most skew products over Anosov systems have the accessibility
property, and also Brin with Gromov and Karcher \cite{brgr},
\cite{brka} proved the accessibility property for some frame flows.
But by that time none of those works guaranteed that this
accessibility was stable. The first case where the accessibility
property was guaranteed to be stable was for the geodesic flows on
surfaces of constant negative curvature, \cite{grpush}, after this
comes the variable curvature, \cite{wi}, the contact Anosov flows,
\cite{kako}, and then in \cite{push2} it was proven that if $E^s$
and $E^u$ are $C^1$ then accessibility implies stable accessibility
which implies that lots of examples are stably accessible, for
example, all the affine diffeomorphisms having the accessibility
property. Then conjecture \ref{accesisgeneric} was proven in the
context of skew products, \cite{buwi1}, then in the context of
Anosov flows by \cite{bupuwi}, using a result of \cite{pl}. In
\cite{shwi1} it was proven that some examples can be approximated by
stable accessible ones. Then in \cite{nito} it was proven when $\dim
E^c=1$ and some more technical requirements that involve the
existence of two nearby closed invariant central leaves and that any
two central leaves can be joined by an $su$-path. Didier, in
\cite{di} proved that accessibility implies stable accessibility if
$\dim(E^c)=1$. On the other hand, in \cite{frh} it is proven that
some ergodic linear automorphisms on tori are stably essentially
accessible, but as they are not accessible then they will be not
stably accessible.


\begin{subsection}{The differentiable case.}
In the case when $E^s$ and $E^u$ are differentiable the behavior of
the accessibility classes are well understood. Essentially the
accessibility classes form a stratification of the manifold. Indeed
the following is true:
\begin{theorem}\cite{su1}, \cite{st}, \cite{push2}
If $E$ and $F$ are differentiable then the accessibility classes
have a natural smooth manifold structure such that the inclusion is
a differentiable map. Moreover the accessibility classes and their
dimension vary lower semi-continuously w.r.t. points in the
manifold.
\end{theorem}
The semi-continuous variation we are talking about above is that if
$x_n\to x$ then the limit of the accessibility classes through $x_n$
contains the accessibility class of $x$ and the limit of the
dimension of the accessibility classes of $x_n$ is larger than the
dimension of the accessibility class of $x$.
\begin{problem}
Prove that if $E^s\oplus E^u$ is differentiable then the
accessibility classes behave as when $E^s$ and $E^u$ are
differentiable.
\end{problem}
\begin{problem}
Prove that, in the general case, the accessibility classes are
topological manifolds that vary semi-continuously as well as their
dimension. Prove that, with bunching, they are indeed smooth
manifolds.
\end{problem}

When the central dimension is $1$ we give a solution to this problem
in \cite{rhrhur1}, see subsection \ref{accescentral1d}. When the
central direction is two dimensional we think that a proof in the
lines of \cite{frh} should be possible. For the general case we
think that the approach in \cite{frh} should be useful, but more
machinery would be needed.\par
More generally, we say that a topological space is {\de
topologically locally homogeneous} if given two points, there is a
local homeomorphism from a neighborhood of one onto a neighborhood
of the other. There are lots of unsolved problems on this subject,
see for example \cite{mo}, \cite{mozi}, \cite{huwa}, but we want to
address the following:
\begin{problem}
Let $X$ be a subset of $\R^n$ and assume it is topologically locally
homogeneous. Assume also that the local homeomorphisms extend to
local diffeomorphisms of $\R^n$ and that they are diffeotopic to the
identity through local diffeomorphisms preserving $X$. Prove that
$X$ is indeed a differentiable manifold.
\end{problem}

With respect to perturbations we have the following,
\begin{theorem}\cite{gra}, \cite{push2}\label{accessopenc1bund}
Let $E$ and $F$ be integrable $C^1$ bundles with the accessibility
property. If $E'$ and $F'$ are uniquely integrable $C^0$ bundles,
$C^0$ close to $E$ and $F$ then $(E',F')$ have the accessibility
property.
\end{theorem}
Thus, accessibility is stable in this setting.
Moreover, having the accessibility property in control theory is a
generic property,
\begin{theorem}\label{accessgencontrol}\cite{lo}
The $C^r$ pairs $(E,F)$, $r\geq 1$, of uniquely integrable plane
fields that have the accessibility property form an open and dense
set.
\end{theorem}
Conjecture \ref{accesisgeneric} looks for the analogous of theorem
\ref{accessgencontrol} in the setting of partially hyperbolic
systems. Observe that the main difficulty in the partially
hyperbolic setting is that we perturb $f$ instead of the bundles
themselves. Maybe the following is not that hard to prove:
\begin{problem}\label{accessaproxprob}
If $E^s\oplus E^u$ is smooth then $f$ can be approximated by a
stably accessible diffeomorphism. Or maybe assuming that $E^s$ and
$E^u$ are smooth.
\end{problem}
Related to this problem we can mention the result of Shub and
Wilkinson \cite{shwi1} that solves problem \ref{accessaproxprob} in
two types of cases. In both cases it is required the integrability
of $E^s\oplus E^u$ and some kind of global product structure of the
center manifolds times the $su$-manifolds.

\begin{subsubsection}{$C^1$ density of stable accessibility}
As we have already said Conjecture \ref{accesisgeneric} was
completely solved in the $C^1$ category by D. Dolgopyat and A.
Wilkinson.

\begin{theorem}\label{dolgowilkinson}\cite{dowi} Stable accessibility is $C^1$ dense among the
partially hyperbolic diffeomorphisms, volume preserving or not.
\end{theorem}

Their strategy is the following: Let $f$ be a partially hyperbolic
diffeomorphism. First of all, they find a finite collection of small
disks, $\{D_j\}$, approximately tangent to the center direction in
such a way that $f$ is accessible modulo these disks. More
precisely, for $x,y\in M$ there is a finite sequence of $su$-paths
such that the first path begins at $x$, the last path ends at $y$
and all the other starting and ending points of the sequence of
paths belong to a disk in $\{D_j\}$, each path beginning in the disk
where the previous path ends. It is not difficult to show, perhaps
enlarging the size of the disks a little bit, that accessibility
modulo $\{D_j\}$ is an open property. After that they perturb $f$ in
a small neighborhood of $\cup D_j$ in such a way that they obtain
(stably) that any pair of points of any disk $D_j$ of the collection
can be joined by an $su$-path. These local accessibility together
the accessibility modulo $\{D_j\}$ give the accessibility for $f$.
To prove the local accessibility in a disk some control on the
effect of the  perturbation of $f$ on the strong bundles is needed.
They show that the contribution of a $C^1$ perturbation in a small
neighborhood of a point $x$ to $E^s(x)$ and $E^u(x)$ is larger than
the contribution of the rest of the orbit. This seems to be false
for a $C^2$ perturbation.
\end{subsubsection}

\end{subsection}

\begin{subsection}{The case $\dim E^c=1$.}\label{accescentral1d}
In \cite{rhrhur1} it was proven conjecture \ref{accesisgeneric} for
the volume preserving case when $\dim E^c=1$. In \cite{burhrhtaur}
the volume preserving condition is removed. Let us see how it works.

The idea is first to understand how do the accessibility classes
behave. In fact the first lemma is very simple, completely general
and it does not need the $\dim E^c=1$ assumption.
\begin{lema}
The following properties are equivalent:
\begin{enumerate}[i)]
\item $C(x)$ is open,
\item $C(x)$ has nonempty interior,
\item the intersection of $C(x)$ with some $c$-dimensional manifold transversal to $E^s\oplus E^u$ has nonempty interior.
\end{enumerate}
\end{lema}

Then, we define the sets $U(f)=\{x;\, C(x)\, \mbox{is open}\}$ and
$\Gamma(f)$ the complement of $U(f)$. $U(f)$ is clearly open and
invariant and hence $\Gamma(f)$ is closed and invariant. By
connectedness of $M$, accessibility means that $U(f)=M$. When a
point is in $\Gamma(f)$ we say that $E^s\oplus E^u$ is integrable.
This is motivated by the following:
\begin{theorem}
$\Gamma(f)$ is laminated by the accessibility classes, that is, the
partition by accessibility classes restricted to $\Gamma(f)$ form a
lamination. Moreover the map $f\to\Gamma(f)$ varies upper
semi-continuously w.r.t. partially hyperbolic $f$ , i.e. if
$f_n\to_{C^1} f$ and $f$ is partially hyperbolic then
$\limsup\Gamma(f_n)\subset\Gamma(f)$.
\end{theorem}
Observe that the semi-continuity of $\Gamma$ automatically
guarantees that the accessibility property is stable, see also
\cite{di}.

So that the accessibility classes are extremely well behaved. We
think that an analogous result should be true in quite full
generality, for example something like:
\begin{problem}
Define $\Gamma_k(f)$, $k=u+s,\dots,N$, to be the set of points $X$
where the dimension of $C(x)$ is $k$ and
$\Lambda_k(f)=\bigcup_{j=u+s}^k\Gamma_k(f)$. Then, $\Lambda_k(f)$ is
closed for every $k$, $\Gamma_j$ is laminated, for every $j$ and
they fit to build a stratification on $\Lambda_k$. Moreover,
$f\to\Lambda_k(f)$ varies semi-continuously.
\end{problem}
The next stage is to know how to break the integrability. To this
end we prove the following:
\begin{proposition}
Let $\Res$ be the set of partially hyperbolic diffeomorphism such
that $Per(f)\subset U(f)$, that is the diffeomorphism whose
accessibility classes at periodic points are open. Then $\Res$ is
$C^{\infty}$ dense.
\end{proposition}
The proof follows the lines of Kupka-Smale theorem once we know how
to make a perturbation at a periodic point to open its accessibility
class. To make the perturbation smooth, the idea is to find where to
make the perturbation. So we have the following lemma that is quite
general and that guarantees that there are lots of points that do
not return to themselves. Let us denote with $\F_{\epsilon}(x)$ the
$\epsilon$ ball of $x$ in the leaf $\F(x)$. Given $x\in M$ let
$\Sigma$ be a transversal to $\F$ at some fixed angle with $\F$. Let
us define
$A_{\epsilon}(x)=\F_{4\epsilon}(\Sigma)\setminus\F_{\epsilon}(\Sigma)$.

\begin{lema}[Keepaway lemma]\label{nonreturning}
Let $f$ be a diffeomorphism preserving a foliation $\F$. Assume that
$m(Df|\F)\geq\lambda>1$ and let $N>0$ be such that $\lambda^N\geq
4$. Then for every $\epsilon>0$ and every $x\in M$, if
$f^j(A_{\epsilon}(x))\cap B_{\epsilon}(x)=\emptyset$ for $j=1,\dots
N$, then there is $y\in\F_{\epsilon}(x)$ such that $d(f^n(y),x)\geq
\epsilon$ for every $n\geq 0$.
\end{lema}
The idea is then to take a periodic point $p$, and a non-returning
point $z\in W^u(p)$ close to $p$. If $p$ is in $\Gamma(f)$ we can
take a three legged $su$-path beginning at $z$ and ending at $p$. We
can suitably choose it in such a way that the breaking points of the
$su$-path are also non-returning, for the past and/or for the
future. Then, essentially any push supported in a small ball around
$z$, transversal to the $E^s\oplus E^u(z)$ direction, will not
change the three legged $su$-path, but will change the unstable leaf
at $z$ and hence it will break the integrability at $p$.\par
So the idea is to find periodic points to prove the following:
\begin{proposition}\label{r=aub}
$\Res=\A\cup\B$ where $\A$ is the set of diffeomorphism having the
accessibility property and $\B$ is the set of diffeomorphism with
$E^s\oplus E^u$ integrable, and having no periodic points, i.e.
$\Gamma(f)=M$ and $Per(f)=\emptyset$.
\end{proposition}
The end of the proof of the denseness of the accessibility property
is just to notice that the set $\B$ is closed by the semi-continuity
of $\Gamma(f)$ and is nowhere dense by the same type of perturbation
carried on periodic points, but now on any other point.\par
To prove proposition \ref{r=aub} it is used the following
proposition which will be quite useful in the description of the
accessibility partition that is carried out in \cite{rhrhur3}, see
also next subsection. Given a compact invariant, $su$-saturated set
$K\subset\Gamma(f)$ let us define the central boundary of $K$,
$\partial^cK$ as the set of points $x\in K$ such that for any
central curve $W$, $x$ is in the boundary of an interval of $K^c\cap
W$.
\begin{proposition}\label{abundper}
Let us assume that $\Omega(f)=M$. If $K\subset\Gamma(f)$ is a
compact invariant, $su$-saturated set then $\partial^cK$ is
$su$-saturated and the set of periodic points in $\partial^cK$ is
dense in $\partial^cK$.
\end{proposition}

The proof of this proposition uses heavily that $\Omega(f)=M$. The
case $\Omega(f)\neq M$ needs another type of argument, exploiting in
a  more subtle way the semi-continuity of $f\to\Gamma(f)$, see
\cite{burhrhtaur} for more details.

\end{subsection}

\begin{subsection}{Some special cases.}
In this subsection we shall exploit proposition \ref{abundper} and
find some interesting description of the partition into
accessibility classes when the unstable (or the stable) manifold has
also dimension $1$.

The main theorem here is an application of proposition
\ref{abundper} and Franks work on codimension one Anosov
diffeomorphisms. To this end let us give first some example. Let $A$
be a codimension one Anosov linear automorphism on $\T^N$ and let
$B$ be another linear automorphism commuting with $A$. Let $M_B$ be
the manifold that is the quotient $\T^N\times\R/\sim$ where
$(x,t)\sim (y,s)$ if and only if $B^nx=y$ and $t=s+n$. Then on $M_B$
can be defined lots of partially hyperbolic systems related to $A$,
for instance, any diffeomorphism $F:\T^N\times [0,1]\to \T^N\times
[0,1]$ such that $F|\T^N\times\{0\}=F|\T^N\times\{1\}$ is homotopic
to $A$ defines one of them. Then we have the following:

\begin{theorem}\label{tori}
Let $F:M\to M$ be a partially hyperbolic diffeomorphism on a compact
manifold $M$. Assume that $\Omega(F)=M$, $\dim E^u=1$ and $\dim
E^c=1$. Then either $M=M_B$ and $F$ is as one of the examples above,
or $E^s\oplus E^u$ is integrable, or $F$ has the accessibility
property, or $M$ decompose as the union $M=\bigcup_{i=1}^n M_i$
where $M_i$ is a compact manifold whose boundary is a finite union
of tori, each of which injects in homotopy, and $M_i\cap
M_j\subset\partial M_i\cap\partial M_j$ for $i\neq j$. $M_i$ is
$su$-saturated, $F^{k_i}$ invariant for some $k_i>0$ and its
boundary are accessibility classes. For each $i$, either the
interior of $M_i$ is itself an accessibility class or $M_i$ is
homeomorphic to a torus times an interval. Finally, if for each $i$
we take the $M_i$ that are homeomorphic to torus times interval in a
maximal fashion then $n$, the number of elements in the
decomposition is less than the first Betti number of $M$, i.e.
$n\leq\dim H_1(M,\R)$.
\end{theorem}
When $M$ is foliated by the accessibility classes, it can be shown
that this $su$ foliation has the following minimality property, the
unique nonempty open $F-$invariant and $su-$saturated set is the
whole $M$.

Notice that in order to have a nontrivial decomposition one needs to
have at least more than one invariant torus injecting in $\pi_1(M)$
and with the dynamics on it being Anosov. This is a highly
nontrivial topological restriction. On the other hand, if all the
components are of the form torus times interval then the dynamics of
$F$ restricted to any of the boundary components are all conjugated
to the same linear map $A$ and $M=M_B$ as one of the examples above.

Finally, observe that when the decomposition is trivial, then either
the system has the accessibility property or $E^s\oplus E^u$
integrates to a foliation that jointly with $F$ have a minimal
property. Aside from some linear Anosov diffeomorphisms on torus, we
do not know of any example of this last case. In the volume
preserving case, we believe that this last case is already enough to
guaranty ergodicity. We want to put the following related problem:

\begin{problem}
Given a $C^2$ codimension one minimal foliations on a compact
manifold $M$, there is essentially a finite number of measurable
saturated subsets of $M$, i.e. the $\sigma-$algebra of measurable
saturated sets is finite mod $0$ w.r.t. Lebesgue.
\end{problem}

We thought that in fact the foliation should be ergodic w.r.t.
Lebesgue measure, but it seems that some construction beginning with
some minimal but non-ergodic flow on a surface could be carried out
to build some non-ergodic minimal foliation. In any case, we think
that  in the particular case that the codimension one foliation
comes from this partially hyperbolic setting, then the
$\sigma$-algebra should be trivial. But it is still not clear in
this case why the $su$-foliation is $C^2$, see problem
\ref{difofsuhol}.

We finish with another problems:
\begin{problem}
Prove that on most three dimensional manifolds, partial
hyperbolicity already implies ergodicity.
\end{problem}

\begin{problem}
If $f:M\to M$ is an Anosov diffeomorphism on a complete riemannian
manifold $M$. Is it true that if $\Omega(f)=M$ then $M$ is compact?
\end{problem}
When $f$ is a codimension one Anosov diffeomorphism the answer is
yes and that is one of the ingredients to the proof of theorem
\ref{tori}.
\end{subsection}
\begin{subsubsection}{Three dimensional nil-manifolds}\label{sacksteder}
Here we shall see what happens in one of the first algebraic,
non-trivial examples. It appears in Sacksteder work \cite{sa} where
he proves its ergodicity using the accessibility property.

Let $\He$ be the Heisenberg group of upper triangular $3\times 3$
matrices with ones in the diagonal. This is the non-abelian
nilpotent simply connected three dimensional Lie group. We may
identify $\He$ with the pairs $(\x,y)$ where $\x=(x_1,x_2)\in\R^2$,
$y\in\R$, $(\x,y)\cdot({\mathbf a},b)=(\x+a,y+b+x_1a_2)$ and
$(\x,y)^{-1}=(-\x,x_1x_2-y)$. For $(\x,y)$ and $({\mathbf a},b)$ in
$\He$, their commutator is $[(\x,y),({\mathbf
a},b)]=(0,x_1a_2-a_1x_2)$. Hence $(\x,y)$ commutes with $({\mathbf
a},b)$ if and only if $\x$ and ${\mathbf a}$ are colinear. We have
also the projection $p:\He\to\R^2$, $p(\x,y)=\x$ which is also an
homomorphism.

If we denote with $\he$ the Lie algebra of $\He$, then $\he$
corresponds with the upper triangular matrices with zeros in the
diagonal. We may also identify $\he$ with the pairs $(\x,y)$ where
$\x=(x_1,x_2)\in\R^2$, $y\in\R$. We have the exponential map
$\exp:\he\to\He$ given by $\exp(\x,y)=(\x,y+\frac{1}{2}x_1x_2)$,
$\exp$ is one to one and onto. Its inverse, the logarithm,
$\log:\He\to\he$ is given by $\log(\x,y)=(\x,y-\frac{1}{2}x_1x_2)$
and the Lie bracket is given also by $[(\x,y),({\mathbf
a},b)]=(0,x_1a_2-a_1x_2)$.

The homomorphisms from $\He$ to $\He$ are of the form
$L(\x,y)=\left(A\x,l(\x,y)\right)$, where
$$
A=\left(\begin{array}{cc}
a & b\\
c & d
\end{array}\right)\;\;\;\;\mbox{and}\;\;\;\; l(\x,y)=\alpha x_1+\beta x_2+\det(A)y+\frac{ac}{2}x_1^2+\frac{bd}{2}x_2^2+bcx_1x_2.
$$

If we denote with $\hat L:\he\to\he$, $\hat L=D_0 L$, it is induced
by the matrix
\begin{displaymath}
\hat L=\left(
\begin{array}{ccc}
a & b & 0\\
c & d & 0\\
\alpha &\beta & \det(A)
\end{array}\right)=\left(
\begin{array}{cc}
A         & 0\\
{\mathbf v} & \det(A)
\end{array}\right)
\end{displaymath}
where ${\mathbf v}=(\alpha,\beta)$ and it follows that
$\exp\left(\hat L(\x,y)\right)=L\left(\exp(\x,y)\right)$.

The centralizer of $\He$ coincides with its first commutator, i.e.
$Z(\He)=[\He,\He]=\He_1$ which consists of the elements of the form
$(0,y)$. Any homomorphism from $\He$ to $\He$ must leave $\He_1$
invariant. Similarly $\he_1=[\he,\he]$ also consists of the elements
of the form $(0,y)$. The automorphisms of $\He$ are exactly the ones
with $\det(A)\neq 0$.

Any lattice in $\He$ is isomorphic to
$\Gamma_k=\{(\x,y):\x\in\Z^2,y\in \frac{1}{k}\Z\}$, for $k$ a
positive integer and the automorphisms leaving $\Gamma_k$ invariant
are the ones with $A\in GL(2,\Z)$ (the matrices with integral
entries and determinant $\pm 1$) and $\alpha,\beta\in\frac{1}{k}\Z$.
On the other hand, every automorphisms of $\Gamma_k$ extends to an
automorphism of $\He$.
\begin{lema}\label{centralizer}
If $T$ is a subgroup of $\Gamma_k$ isomorphic to $\Z^2$, $T\cap \He_1\neq\{(0,0)\}$.
\end{lema}
We define the quotient compact nil-manifold $N_k=\He/\Gamma_k$ by
the relation $(\x,y)\sim (a,b)$ if and only if
$(\x,y)^{-1}\cdot({\mathbf a},b)\in\Gamma_k$. The first homotopy
group is $\pi_1(N_k)=\Gamma_k$ and two maps of $N_k$ to itself are
homotopic if and only if their action on $\Gamma_k$ coincide.
Moreover, any map from $N_k$ to itself is homotopic to an
automorphism as the described above leaving $\Gamma_k$ invariant.
The projection $p:\He\to\R^2$ descends to a projection
$p:N_k\to\T^2$ and also $p$ serves as semiconjugacy between any
automorphism $L$ and its corresponding matrix $A\in GL(2,\Z)$.\par
Notice that lemma \ref{centralizer} plus theorem \ref{tori} and its
comments imply that any partially hyperbolic diffeomorphism on $N_k$
either has the accessibility property or $E^s\oplus E^u$ integrates
to a minimal foliation. We think that this second possibility does
not exist.\par
Given an automorphism $L$ of $N_k$, it is partially hyperbolic if
and only if the associated matrix $A\in GL(2,\Z)$ is hyperbolic. In
this case, by taking some finite covering it is possibly to make
${\mathbf v}=0$, but we will not use this fact here.

Let us denote the invariant subspaces for $A$ with $E^s_A,
E^u_A\subset\R^2$ with corresponding eigenvalues $\lambda^s$ and
$\lambda^u$. Then the invariant subspaces for $D_0L=\hat L$ are
$$
E^s_L=\left\{\left({\mathbf w},\frac{{\mathbf v}\cdot{\mathbf
w}}{\lambda^s-\det(A)}\right)\;:\;{\mathbf w}\in E^s_A\right\},
$$
and similarly
$$
E^u_L=\left\{\left({\mathbf w},\frac{{\mathbf v}\cdot{\mathbf
w}}{\lambda^u-\det(A)}\right)\;:\;{\mathbf w}\in E^u_A\right\}.
$$
Finally, $E^c_L=\he_1$ and the invariant bundles are formed by the
translates of the invariant spaces at $0$. Similarly, the translates
of $E^{\sigma}_L$, $\sigma=s, u, c$ project onto the invariant
foliations in $N_k$. Observe that $E^c_L$ projects onto a circle and
hence the projection of the translates of $E^c_L$ is a foliation by
circles. Moreover, these circles are collapsed by the projection $p$
and hence the central foliation is a nontrivial fibration with base
$\T^2$ and fiber $\S^1$.

Finally,  the accessibility class of $0$  will be the smallest
subgroup containing $E^s_L$, $E^u_L$ and being closed under the
bracket operation. Thus, its lift to the Lie algebra will be the
smallest Lie sub-algebra containing $E^s_L$ and $E^u_L$ that also
must contain $\he_1=[E^s_L,E^u_L]$ and hence equals $\he$. So that
any such $L$ has the accessibility property.

\end{subsubsection}
\begin{subsubsection}{The general affine case}
Let $f:G/B\to G/B$ be  a partially hyperbolic affine diffeomorphism.
Recall, subsection \ref{afdif}, that this is equivalent to say that
$\he\not\subset\go{b}$ where $\he$ is the hyperbolic subalgebra of
$f$. Let $\bar{f}:G\to G$, $\bar{f}=L_g\circ A$ be the affine
diffeomorphism covering $f$ and $\au(f):\g\to\g$ the corresponding
automorphism of its Lie algebra. Observe that the partition into
accessibility classes for $\bar{f}$ corresponds exactly to the
translates of $H$, the hyperbolic subgroup of $\bar{f}$ and hence it
is a foliation, moreover it comes from the left action of $H$ on
$G$. The following is also true:
\begin{theorem}\cite{bresh}, \cite{push4}
If $f:G/B\to G/B$ is a partially hyperbolic affine diffeomorphism
then its accessibility classes are the orbits of the left action of
$H$ on $G/B$. Hence $f$ has the accessibility property if and only
if $\he+\go{b}=\g$. Moreover, $f$ has the essential accessibility
property if and only if $\overline{HB}=G$.
\end{theorem}
With respect to the way in which the accessibility classes behave
after perturbations, the following is proven by Starkov in the
appendix of \cite{push4},
\begin{proposition}\label{pertlefttrans}
There exists a neighborhood $O(g)\subset G$ such that $H$ is
contained in the hyperbolic subgroup $H_x$ of $f_x=L_x\circ A$ for
every $x\in O(g)$.
\end{proposition}
So that after perturbation by left translation the accessibility
classes can only grow. It would be interesting to know in which
cases it is possible to get accessibility, or essential
accessibility only by applying left translation, i.e.
\begin{problem}
Given an automorphism of a finite volume homogeneous space,
$A:G/B\to G/B$ let us define $G_i\subset G$, $i=\dim G^s+\dim
G^u+\dim B, \dots, d=\dim G$ by $G_i=\{g\in G\;\mbox{such that}\;
\dim \overline{H_gB}=i\}$ where $H_g$ is the hyperbolic subgroup of
$L_g\circ A$. This gives a partition of $G$. How is the structure of
this partition?, How do $G_i$ look like, are they manifolds?, For
which automorphisms is $G_d\neq\emptyset$? For which automorphisms
is $G_i\neq\emptyset$ for every $i$? Is this possible? What is their
relation with the algebraic properties of $G$? Analyze also the
analogous subsets but for $\he_g+\go{b}$ instead of
$\overline{H_gB}$, do they coincide?
\end{problem}
Hence in the affine case, the accessibility classes behave the best
possible, also under perturbations. As we said, we think that this
should be true in the general context.
\begin{problem}
Study the affine diffeomorphism $f$ on $G/B$ having $\he+\go{b}$ of
codimension one and their perturbations. It seems that some
description of all possible cases should be plausible, at least with
low dimensional center space. What about codimension $2$? What about
one and two dimensional center space?
\end{problem}
\end{subsubsection}

\begin{subsection}{Stable ergodicity of toral automorphisms}\label{ergodicautomorphism}
In \cite{frh} it is given a partial answer to problem
\ref{ergodicautomorphismproblem}. It is proven the stable ergodicity
for some toral automorphism with two dimensional center bundle and
some extra assumption. This extra assumption essentially asks for
the irreducibility of the matrix $A$ defining the linear
automorphisms. In fact, the assumption is that all the powers of $A$
have irreducible characteristic polynomial. This irreducibility
condition is not very much restrictive because every toral
automorphisms can be essentially decomposed by blocks satisfying
this irreducibility condition. Thus, if one understands what the
behavior is for this type of linear automorphisms and its
perturbations, its seems likely that one can go through to the
general case. The restrictive assumption here is the two
dimensionality of the center bundle. Let us formulate the results:

\begin{theorem}\label{n6}
Every linear automorphism of $\T^N$ with the hypothesis listed above
is $C^5$-stably ergodic if $N\geq 6$.
\end{theorem}
When $N=4$ we have the following:
\begin{theorem}\label{n4}
Every linear automorphisms of $\T^4$ with the hypothesis listed
above is $C^{22}$-stably ergodic.
\end{theorem}
Here $C^r$-stable ergodicity means that the perturbations are made
in the $C^r$ topology.\par
Moreover, all ergodic linear automorphisms acting on $\T^4$ are
either Anosov or satisfiy the extra assumptions; and all ergodic
linear automorphisms acting on $\T^2$, $\T^3$ and $\T^5$ are Anosov.
Hence, we have the obvious corollary:
\begin{theorem}\label{sol}
Every ergodic linear automorphism of $\T^N$ is stably ergodic for
$N\leq 5$.
\end{theorem}
Before giving a rough idea of how the proof goes, let us spend a few
words about the high differentiability used. In fact most of the
proof follows only with a $C^1$ assumption, but at some places it is
used KAM linearizing theorems that make use of high
differentiability assumption. When $N\geq 6$ it is used the standard
linearizing theorem of Arnold and Moser but for $N=4$ things look a
little bit different so it is needed to adapt a linearizing theorem
of Moser on commuting diffeomorphisms of the circle to this setting
but then the differentiability increases to $22$. We think that the
differentiability could be improved without new techniques to allow
$C^1$ perturbations but by $C^r$ diffeomorphism.

\begin{problem}
Are the ergodic automorphisms of tori stably Bernoulli? Maybe a
little simpler, in dimension 4 when the $E^s\oplus E^u$ bundle is
integrable, is $f$ Bernoulli?
\end{problem}

Let us go now into the proof. As one can imagine, the idea is to
prove the essential accessibility of the perturbations. Observe that
for the linear map, the bundle $E^s\oplus E^u$ integrates to a
foliation by planes that coincide with the partition by
accessibility classes. So there is no hope of getting stable
accessibility, and one needs to find a way of distinguishing when
there is accessibility.\par
The first step is to prove that the perturbation, in some sense
looks much like the linear case. To this end, it is much more useful
to work in the universal covering $\R^N$. One of the most useful
properties in the linear case is that all the invariant foliations
are foliations by planes. So any two leaves of the same foliation
are parallel and two leaves of different foliations always
intersect. For the perturbation we have that all the distributions
but the $su$ are integrable. Moreover the center, center-stable and
center-unstable leaves for $f$ stay in a bounded tubular
neighborhood of the corresponding for the linear. But for the strong
foliations this is no longer the case, typically the leaves  will
not stay in a bounded tubular neighborhood of the linear ones.
Nevertheless, if one works with logarithmic type tubular
neighborhoods, then one gets that the leaves of the strong
foliations stay in that logarithmic type tubular neighborhoods of
the linear ones. This is enough to guarantee that any leaf of the
stable foliation intersects any leaf of the center-unstable at
exactly one point, and the same holds for the unstable foliation.
Moreover, this also allows us to define the asymptotic direction for
the strong foliations and to see that they coincide with the linear
case. So one can define global stable holonomies between two
center-unstable leaves and global unstable holonomies between two
center-stable leaves.\par
The second step is to study the stable and unstable holonomies. It
turns out that when one restricts the stable holonomy to center
leaves (whenever it makes sense) then they are differentiable.
Moreover, one sees that if the center leaves are not much far away
then the stable holonomy between center leaves are close to the ones
of the linear in the $C^r$ uniform topology in the whole center leaf
if the perturbation is $C^r$ small. Besides, when the center leaves
are far, one can still measure the Lipschitz constant of the stable
holonomy and see that is no worse than a small power of the distance
between the leaves. With this control on the growth of the Lipschitz
constant, one can measure some type of growth of the volume of
accessibility classes. For example, if one takes a central ball of
radius $\epsilon$ and then take all the unstable balls of radius
$1/\epsilon^{\beta}$ and on this set, the stable balls of radius
$\epsilon$,  it follows that this set is open, and for some choice
of $\beta$, its volume is more or less $1/\epsilon^{\gamma}$ for
some $\gamma>0$. This gives some type of recurrence for the
accessibility classes.\par
%
In the third step it is shown, using the listed properties and that
the central direction is two dimensional, that the partition by
accessibility classes is essentially minimal. In fact it is proven
that the only open, $f$-invariant and $su$-saturated sets are the
whole manifold and the empty set. The idea of the proof is to
recover same flavor of the Denjoy argument for the rotations on the
circle. The information on the growth of the volume will take the
place, in some sense, of the Denjoy-Koksma inequality. One also make
use of the diophantine property of the asymptotic directions of the
strong foliations (which are the same as the linear one). Finally,
the criterion used to prove that such a non-empty open set is the
whole manifold is to have the same homotopy type of the torus, and
as the torus is a $K(\Z^N,1)$ Eilenberg-MacLane space this is done
by proving that all the homotopy groups $\pi_k$, $k\geq 2$ are
trivial. Here is, maybe, the crucial step where the two
dimensionality of the center direction appears. Indeed, to make all
the homotopy groups trivial, it is seen that the homotopy groups of
an $su$-saturated set are the same as the ones of its intersection
with a central leaf and hence, as the $\pi_k$ of any subset of the
plane is trivial if $k\geq 2$, one only has to deal with the $\pi_0$
and $\pi_1$. Observe that this fact is no longer true in higher
dimensions ($\pi_3(S^2)\neq 0$). The $\pi_0$ is treated easily, it
is the $\pi_1$ that consumes the biggest effort.\par
Observe that by the above minimality property it is quite easy to
make the accessibility property appear, for example, if there is an
accessibility class with nonempty interior then, as the system is
volume preserving, it is not hard to see that it should be open and
essentially invariant and hence we get the accessibility. So, to get
accessibility it is enough to get some open accessibility class
somewhere. But as we have said, there is no hope to get always
accessibility, so it is needed to know what happens when there is no
accessibility, and that is step 4.\par
In the fourth step, it is shown that the accessibility classes are
(topological) manifolds and moreover, their dimension vary
semi-continuously. Although at this stage it is used the fact that
the central dimension is two, we believe that this is true in full
generality. Here one mostly works with the intersection of an
accessibility class with a center leaf which we call {\de central
accessibility class} and call the dimension of this intersection the
{\de central dimension}. We get for example that the set of points
whose accessibility class has zero central dimension is a closed,
$f$-invariant, $su$-saturated set and hence, by the minimality, is
either the whole manifold or the empty set. If it were empty then
either would we have that there is an accessibility class that has
central dimension two and hence should be open and hence we get the
accessibility property, or else all accessibility classes would have
central dimension one and hence the central accessibility classes
would be curves, in particular this would be the case for a fix
point (there is always a fixed point) that needs to have complex
eigenvalues if the perturbation is small. So we would have an
invariant curve through a fixed point with a complex eigenvalue,
thus this curve should spiral. But then we know that the
accessibility classes are homogeneous (a neighborhood of each point
is homeomorphic to some neighborhood of any other point) and
moreover the homeomorphisms that make it homogeneous are holonomies
between central leaves and hence they are diffeomorphisms. So we get
that the curve should spiral at all its points and this is
impossible. Thus, we get that if we do not have the accessibility
property then the central accessibility classes are zero
dimensional, that is, they are points. Recall that we are working in
the universal covering, thus if we fix a central leaf, we get that
the intersection of an accessibility class with this central leaf is
just a point. So we finally get what we were looking for, what
happens when we do not have the accessibility property.\par
In the final step we have to deal with the case where we do not have
the accessibility property and we have to prove the essential
accessibility property. At this stage the two dimensionality of the
central foliation is not needed. The step above allows us to define
the $su$-holonomy between center leaves simply as the intersection
of the accessibility class with the leaf. By the properties we got
in the second step, this holonomies are $C^r$ and $C^r$-close to the
linear ones. This property and the fact that the asymptotic
directions are diophantine allow us to use the KAM scheme to get a
smooth conjugacy between the perturbed partition by accessibility
classes and the linear one, thus getting the essential accessibility
for the perturbed system and hence the stable ergodicity. The way we
use the KAM scheme is by taking a global transversal to the
partition by accessibility classes and hence to get a $\Z^N$ acting
on it. Thus we take $2$ elements of the action and build a
two-dimensional torus on it and we still have an action of
$\Z^{N-2}$ on this two dimensional torus. At this point, the case of
dimension $4$ and dimension bigger than $4$ are quite different. In
dimension bigger than four, any element of this action on the torus
has irrational (Diophantine) translation vector with irrational
(Diophantine) slopes and we can find an element close to the
corresponding translation and hence the usual linearizing theorem of
Arnold and Moser applies, see \cite{he}, and hence the
differentiability required is $5$, in fact $4+\alpha$, any
$\alpha>0$ will be enough. But when $N=4$, this is no longer the
case, in fact, already for the linear case any element of the $\Z^2$
action on the two torus have translation vectors that have rational
slopes, we find this an interesting feature that we were not able to
see reflected in the $su-$foliation. In this case we need to use the
$\Z^2$ action in its full strength and we have to prove a theorem on
linearization of commuting diffeomorphism on torus analogous to the
one of Moser in \cite{mos} and to make the algorithm converge we
need the differentiability of at least $21+\alpha$ for some definite
$\alpha>0$ that can be computed. Nevertheless this differentiability
is required  for the convergence of this specific algorithm and it
seems to be far from optimal. Another approach that may improve the
required differentiability is the result of R. Hamilton on the
stability of some foliations.
\end{subsection}
\end{section}

\begin{section}{Ergodic properties of partially hyperbolic systems}\label{ergodicity}
Given a diffeomorphism $g:N\to N$, we say that $g$ is {\de ergodic
with respect to Lebesgue measure} if any invariant set has either
full or null measure; we do not assume a priori that $g$ preserves
Lebesgue measure. More generally, given a partition of a manifold,
we say that {\de the partition is ergodic} w.r.t. Lebesgue measure
if any measurable saturated set has either full or null measure. We
say that a partially hyperbolic diffeomorphism has the {\de
essential accessibility property} if its partition into
accessibility classes is ergodic.

\begin{subsection}{The differentiable case}\cite{ho},\cite{sa},\cite{pushst}
We shall first present a sketch of the proof that some partially
hyperbolic systems are ergodic when the stable and unstable
foliations are smooth. This sketch, popularly known as the Hopf
argument, is essentially the one used in the more general case.
\begin{theorem}\label{ergodicsmooth}
Let $f$ be a partially hyperbolic system and assume that $E^s$ and
$E^u$ are differentiable. Then if $f$ has the essential
accessibility property, $f$ is ergodic.
\end{theorem}
This theorem is a consequence of a more general one that may be
found in \cite{sa}. In \cite{pushst} there is also a proof
for the homogeneous case. This gives essentially the ergodic
decomposition of Lebesgue measure.
\begin{theorem}\label{ergodicdecomposition}
Let $f$ be a diffeomorphism and assume it preserves two smooth
foliations, $\F^s$ and $\F^u$. Assume that vectors tangent to $\F^s$
are exponentially contracted and vectors tangent to $\F^u$ are
exponentially expanded. Then, for each invariant continuous function
$\phi$, $\phi$ is essentially constant along accessibility classes,
that is, there is a full measure set $R$ such that if $x,y\in R$ and
$x$ is in the accessibility class of $y$ then $\phi(x)=\phi(y)$. In
other words, the accessibility classes determine the ergodic
decomposition of Lebesgue measure.
\end{theorem}
The property described in this last theorem is sometimes called the
Mautner's phenomenon.\par
Observe that in this theorem it is not assumed that $f$ is partially
hyperbolic. Nevertheless, it would be interesting to have an example
of a non partially hyperbolic system  where the theorem applies.\par
Instead of giving a proof of theorem \ref{ergodicdecomposition} we
shall go directly to the proof of theorem \ref{ergodicsmooth}.
\begin{proof}
By the Birkhoff ergodic theorem we have that for every $L^1$
function $\phi$, the times or Birkhoff averages
\begin{eqnarray}\label{birkhoff}
\frac{1}{n}\sum_{k=0}^n\phi(f^k(x))
\;\;\;\mbox{and}\;\;\;\frac{1}{n}\sum_{k=0}^n\phi(f^{-k}(x))
\end{eqnarray}
converge a.e., as $n\to +\infty$, to measurable invariant functions
$\tilde{\phi}^+$ and $\tilde{\phi}^-$ respectively. Moreover,
$\tilde{\phi}^-=\tilde{\phi}^+$ a.e. Finally, $f$ is ergodic if and
only if $\tilde{\phi}^+$ is constant a.e. for every continuous
$\phi$. Thus we wish to prove that $\tilde{\phi}^+$ is constant a.e.

Let us call $R^+$ the set where the forward Birkhoff averages are
convergent, and $R^-$ the set where backward Birkhoff averages are
convergent. Observe that by uniform continuity of $\phi$ we have
that $R^+$ is saturated by stable leaves, and $R^-$ is saturated by
unstable leaves. Moreover, if $y$ is in the stable leaf of $x\in
R^+$ then $\tilde{\phi}^+(x)=\tilde{\phi}^+(y)$, similarly if $z$ is
in the unstable leaf of $x\in R^-$ then
$\tilde{\phi}^-(x)=\tilde{\phi}^-(z)$. In other words,
$\tilde{\phi}^+$ is constant along stable leaves and
$\tilde{\phi}^-$ is constant along unstable leaves. Thus, if either
the stable or the unstable foliations where ergodic w.r.t. Lebesgue
measure, then we would get ergodicity of $f$. But a priori, we do
not know if this is the case (see Problem \ref{problemaerg}). Let us
see how we will overcome this difficulty.\par
As Lebesgue measure is invariant, we have that the integrals of
$\tilde{\phi}^+$, $\tilde{\phi}^-$ and $\phi$  are equal and let us
call them $I$. Call $A^+$ the set of $x$ such that
$\tilde{\phi}^+(x)$ is less than $I$, $A^-$ the set of $x$ such that
$\tilde{\phi}^-(x)$ is less than $I$, and observe that $A^+$ and
$A^-$ differ in a null measure set. We need to prove that $A^+$ has
null measure too. So, proving Theorem \ref{ergodicsmooth} is reduced
to proving:
\begin{proposition}\label{densitysat}
The set of Lebesgue density points of an $s$-saturated set, is also
$s$-saturated.
\end{proposition}
This finishes the proof of the theorem because of the following:
Lebesgue density points of sets that differ in a null measure set
are the same, so, if we denote by $D(X)$ the set of Lebesgue density
points of the set $X$, we have that $D(A^+)=D(A^-)$. Now, the
proposition says that $D(A^+)$ is $s$-saturated, and that $D(A^-)$
is $u$-saturated. Hence, $D(A^+)=D(A^-)$ is $s$ and $u$-saturated.
By the essential accessibility property, $D(A^+)$ has either full or
null measure. Hence $A^+$ has either full or null measure, it
clearly cannot have full measure, so we are done.
\end{proof}

Let us go into the proof of proposition \ref{densitysat}
\begin{proof}{of proposition \ref{densitysat}}
The proof of the proposition uses the following general lemma:
\begin{lema}\label{lema.abs.cont}
Let $\F$ be an absolutely continuous foliation and let $\Sigma$ be a
transversal to $\F$. There is a constant $C>0$ such that if $A$ is
an $\F$-saturated set then
$$
\frac{1}{C}m(A)\leq m_{\Sigma}(A\cap\Sigma)\leq Cm(A)
$$
where $m_{\Sigma}$ is Lebesgue measure on $\Sigma$.
\end{lema}
In fact this is essentially the definition of absolute continuity
for a foliation we use. As a consequence of the lemma we have that
$x$ is a Lebesgue density point for a saturated set $A$ if and only
if $x$ is a Lebesgue density point of $A\cap\Sigma$ for the measure
$m_{\Sigma}$. Thus, the proposition follows from this observation
and the fact that the holonomies are $C^1$, and $C^1$
diffeomorphisms preserve Lebesgue density points.
\end{proof}

Let us observe that in the proof of Theorem \ref{ergodicsmooth} the
$C^1$ hypothesis was only used in the proof of proposition
\ref{densitysat}. In fact, this is the proposition that is being
generalized for more general settings in \cite{grpush},
\cite{push2}, \cite{push3}, \cite{buwi2}, \cite{buwi3},
\cite{rhrhur1}. Nevertheless, it is not clear that proposition
\ref{densitysat} remains valid in the general case. In some cases,
instead of using proposition \ref{densitysat}, it is only used the
absolute continuity as in the Anosov case, but in this case
typically it is needed some extra hypothesis in the Lyapunov
exponents, see \cite{budope}; we shall review this item in
subsection \ref{ergodiclyap}.

Observe that the proof of Theorem \ref{ergodicsmooth} works without
change if one only requires, instead of partial hyperbolicity, that
the diffeomorphisms leave invariant two $C^1$ foliations $\F^s$ and
$\F^u$ with the property that points in an $s$-leaf are forward
asymptotic and points in an $s$-leaf are backward asymptotic (not
necessarily exponentially fast).

Let us put some problems related to the proof above.
\begin{problem}\label{problemaerg}
If $f$ is partially hyperbolic and has the essential accessibility
property, are the stable and/or the unstable foliations ergodic
w.r.t. Lebesgue measure? In the homogeneous case it is true,
\cite{pushst}. What about the general $C^1$ case (the stable and
unstable foliations are smooth), or just putting some bunching
condition?
\end{problem}

Let us remark that there is an affirmative answer for the
topological analogous problem, $\dim M=3$: for an open and dense set
of partially hyperbolic systems $f$, if $f$ is robustly transitive,
then the stable or the unstable foliation is minimal \cite{bodiur}.
In \cite{rhrhur2}, the authors extend this result to $\dim M\geq 3$
and $\dim E^c=1$. Also, note that Problem \ref{problemaerg} has an
affirmative answer in the case of Anosov diffeomorphisms.
\begin{problem}\label{ergoddecompprob}
Prove the analogous of theorem \ref{ergodicdecomposition} for
general partially hyperbolic systems.
\end{problem}
When $f$ has some center bunching then problem \ref{ergoddecompprob}
has a positive answer. Let us call $GC(x)=\bigcup_{n\in\Z}C(f^n(x))$
the $f-$saturation of the accessibility classes, then

\begin{theorem}
Let $f$ be a center bunched partially hyperbolic diffeomorphism.
Then the measurable hull of the partition into $GC(x)$ generalized
accessibility classes coincides mod $0$ with the partition into
ergodic components.
\end{theorem}
In other words we will be proving that any $f$-invariant measurable
set coincides mod $0$ with an $f$-invariant $su$-saturated set.
\begin{proof}
By proposition \ref{densitysatbunch} we know that the set of density
points of every essentially $s$-saturated and essentially
$u$-saturated set is $su$-saturated. So we only need to prove that
any invariant set $A$ coincides mod $0$ with an $s$-saturated set
$A^s$ and also with an $u$-saturated set $A^u$. Take a sequence of
continuous functions $\varphi_n$ converging a.e. to $\chi_A$ the
characteristic function of $A$. By Birkhoff ergodic theorem we know
that
$\lim_{N\to+\infty}\frac{1}{N}\sum_{k=0}^{N-1}\varphi_n(f^k(x))\to\tilde{\varphi}^+_n(x)$
for Lebesgue almost every $x$, where $\tilde{\varphi}^+_n$ is a
measurable invariant function. Let as define the sets $ V_n=\{x\in
M\;\mbox{such that}\;\tilde{\varphi}^+_n(x)>\frac{1}{2}\}$. As in
the proof of theorem \ref{ergodicsmooth}, since $\varphi_n$ is
continuous, if $y\in W^s(x)$ then $\tilde{\varphi}^+_n(x)$ exists if
and only if $\tilde{\varphi}^+_n(y)$ exists and in the case they
exist they are equal. Also, by Birkhoff ergodic theorem we know that
$\tilde{\varphi}^+_n\to\tilde{\chi}_A=\chi_A$ a.e. Let
$A^s=\bigcup_{N\geq 0}\bigcap_{n\geq N}V_n$ and observe that it is
$s$-saturated. We claim that $Leb(A\bigtriangleup A^s)=0$. In fact,
let $Z$ be the set of points $x$ such that $\tilde{\varphi}_n(x)$
exists for every $n$ and such that
$\tilde{\varphi}_n(x)\to\chi_A(x)$ and let us see first that $A\cap
Z\subset A^s$. If $x\in A\cap Z$ then we have that
$\tilde{\varphi}_n(x)\to\chi_A(x)$ since $x\in Z$ but as $x$ is also
in $A$ then $\chi_A(x)=1$ and hence there is $N_x$ such that
$\tilde{\varphi}_n(x)>\frac{1}{2}$ for every $n\geq N_x$ hence $x\in
V_n$ for every $n\geq N_x$ and so $x\in A^s$.\par
On the other hand, if $x\in A^s\cap Z$ then we have that
$\tilde{\varphi}_n(x)\to\chi_A(x)$ since $x\in Z$ and as $x$ is also
in $A^s$ there is $N_x$ such that $x\in V_n$ for every $n\geq N_x$.
Hence $\tilde{\varphi}_n(x)>\frac{1}{2}$ for every $n\geq N_x$. So
we get that $\chi_A(x)=\lim\tilde{\varphi}_n(x)\geq\frac{1}{2}$ but
as $\chi_A(x)$ can only be $0$ or $1$ we get that $\chi_A(x)=1$ and
hence $x\in A$. So $A$ coincides mod $0$ with $A^s$ and hence it is
essentially $s$-saturated.\par
The proof that it is essentially $u$-saturated is exactly the same
but putting $N\to -\infty$ in Birkhoff theorem.
\end{proof}
Let us follow with some other problems,
\begin{problem}
Prove ergodicity only assuming that $E^s\oplus E^u$ is
differentiable.
\end{problem}
\begin{problem}
Prove that if  $E^s$ and $E^u$ are differentiable then the
diffeomorphism is approached by an ergodic one, or even by a stably
ergodic one.
\end{problem}
\begin{problem}
Prove that if  $E^s$ and $E^u$ are differentiable and have the
accessibility property then the diffeomorphism is stably ergodic.
\end{problem}


\end{subsection}

\begin{subsection}{Accessibility implies ergodicity}
In this section we shall give an idea of how the proof of the
following theorem goes.
\begin{theorem}\label{bettercenterbunch}\cite{buwi3}
Let $f$ be a $C^2$  partially hyperbolic diffeomorphism, an let all
unit vectors $v^\s\in E^\s$, with $\s=s,c,u$ satisfy the following
inequalities
$$|D_pfv^s|<\nu(p)<\gamma(p)<|D_pfv^c|<\hat\gamma^{-1}(p)<\hat\nu^{-1}(p)<|D_pfv^u|$$
where $\nu,\hat{\nu}<1$. Let $f$ satisfy the following center
bunching conditions \begin{equation}\label{bunching.condition}
\nu<\gamma\hat{\gamma}\qquad\mbox{and}\qquad\hat\nu<\gamma\hat\gamma\end{equation}
Then, if $f$ has the essential accessibility property, it is
ergodic.
\end{theorem}
The main ingredients of the proof of this theorem may be found
already in \cite{grpush}. They were subsequently improved reaching
this final, quite general form. It seems likely that in order to go
further, some new techniques should appear. For example, one of the
main uses of the center bunching condition is to prove
differentiability of the strong foliations when restricted to some
weak ``fake" foliations and the bunching condition is sharp for this
purpose. The proof, a priori, relies heavily on this
differentiability. With this idea we shall put the following
problem:

\begin{problem}
As a step to remove the center bunching assumption, maybe it is
useful to reduce it only to one of the inequalities, for example,
assume only that $\nu<\gamma\hat{\gamma}$. It would be interesting
also to understand what happens in the limit case, for example,
assume that $\nu(p),\hat{\nu}(p)\leq\gamma(p)\hat{\gamma}(p)$ for
every point and equality only holds for a fixed point and at most at
that point the corresponding holonomies are not differentiable. How
should the proof work then? What about assuming only that the stable
foliation is smooth?
\end{problem}
\begin{problem}
Is it true that essential accessibility implies ergodicity only
assuming that the stable and unstable holonomies are $C^1$ when
restricted to center-stables and center-unstables?
\end{problem}


\end{subsection}
\begin{subsection}{A sketch of the proof when $\dim E^c=1$}
Let us give an idea of how we prove, in \cite{rhrhur1}, that
accessibility implies ergodicity when $\dim E^c=1$. As it was said
in p. \pageref{densitysat}, we are reduced to proving the following:
\setcounter{theorem}{2}
\begin{proposition}
The set of Lebesgue density points of an $s$-saturated set, is also
$s$-saturated.
\end{proposition}\setcounter{theorem}{6}
In fact, we shall prove a weaker result, that will be enough for our
purpose. Let us recall that, for $\s=s,u$, an \emph{essentially
$\s$-saturated} set is one that differs from an $\s$-saturated set
only in a set of null measure.
\begin{proposition}\label{densitysatbunch}
The set of Lebesgue density points of an essentially $s-$saturated
and essentially $u-$saturated set is $su-$saturated.
\end{proposition}
That is, Lebesgue density points of essentially $s-$ and essentially
$u-$saturated sets flow through stable and unstable leaves. In
\cite{grpush}, it was suggested that certain shapes called juliennes
would be more natural, rather than merely riemannian balls, in order
to treat preservation of density points. We will follow this line
and use certain solid juliennes instead of balls. Of course, these
new neighborhood bases will define different sets of density points.
Let us say that a point $x$ is a {\de $C_n$-density point} of a set
$X$ if $\{C_n(x)\}_n$ is a local neighborhood basis of $x$, and
$$\lim_{n\to \infty}\frac{m(X\cap C_n(x))}{m(C_n(x))}=1$$
Recall that {\de Lebesgue density points} are $B_{r^n}$-density
point, where $B_{r^n}(x)$ are riemannian balls of radii $0<r^n<1$.
We will be particularly interested in a dynamically defined
neighborhood basis $\{J_n\}_n$, the {\de juliennes}, that consists
of certain local stable and unstable saturations of a small center
arc. For this new neighborhood basis we obtain:
\begin{proposition}
The set of $J_n$-density points of an essentially $s$-saturated set
is $s$-saturated.
\end{proposition}
By changing the neighborhood basis, we have solved the problem of
preserving density points, that is we have established proposition
\ref{densitysat} but for julienne density points. However, we need
to know now what the relationship is between the julienne density
points, and Lebesgue density points. Given a family $\mathcal{M}$ of
measurable sets, let us say that two systems $\{C_n\}_n$ and
$\{E_n\}_n$ are {\de Vitali equivalent} over $\mathcal{M}$, if the
set of $C_n$-density points of $X$ equals (pointwise) the set of
$E_n$-density points of $X$ for all $X\in\mathcal{M}$. That is, if
$D_{C_n}(X)=D_{E_n}(X)$ for all sets $X\in \mathcal{M}$.\par We
obtain
\begin{proposition}\label{jn vitali equivalent lebesgue}
$\{J_n\}$ is Vitali equivalent to Lebesgue over essentially
$u$-saturated sets.
\end{proposition}
Hence, if $X$ is an essentially $s$- and essentially $u$-saturated
set, we have that $D(X)$, the set of Lebesgue density points of $X$,
is an $s$- and $u$-saturated set. As we have said, this is enough to
prove Theorem \ref{ergodicsmooth} for this setting. Indeed, let
$\phi$ be a continuous function, and let $A^+$ the set of points $x$
such that
$$\lim_{n\to+\infty}\frac1n\sum_{k=0}^{n-1}\phi\circ f^k(x)<\int\phi $$
Also let $A^-$ be the set of points such that $\lim_{n\to\infty}
\frac1n\sum_{k=0}^{n-1}\phi\circ f^{-k}(x)<\int\phi$. Then $A^+$ and
$A^-$ differ in a set of null measure, and their density points set
is essentially $s$- and essentially $u$-saturated, so the proof
follows as in the differentiable case to prove that $A^+$ has null
measure, see below Proposition \ref{densitysat}.
\end{subsection}
\begin{subsection}{Juliennes}\label{juliennes}
Let us briefly mention the construction of the juliennes. As one can
infer from the previous paragraph, juliennes are local bases whose
main features are: (1) their density points are equivalent to
Lebesgue over essentially $u$-saturated sets (2) their density
points are $s$-saturated if the set is essentially $s$-saturated.
Let us just try to construct them to fulfill these conditions.\par
As it was said above, juliennes are dynamically defined balls
obtained by locally $su$-saturating center arcs in a certain way.
\par Now, in order that condition (2) be fulfilled, we need that the
julienne density points can flow through stable leaves. If we
restrict ourselves to a center-stable leaf, this is indeed the case,
because the stable holonomy is $C^1$ if the bunching conditions
(\ref{bunching.condition}) are satisfied, and, in particular, if
$\dim E^c=1$. So, if we merely $su$-saturate center arcs of a
certain length $\s^n\in(0,1)$, by arcs of the same length, we would
obtain a ``cubic system" that is easily seen to be equivalent to
Lebesgue, and to preserve density points when restricted to the
center-stable leaf of the center point. This is not enough yet.\par
Indeed, since the global stable holonomy is not $C^1$ in general,
the unstable saturation of the center arc could be very much
distorted by the stable holonomy, thus possibly producing the
appearance of new density points or the disappearance of old ones.
However, if the saturation is dynamically defined, one can bound
this distortion:\par Let us assume for simplicity that $\nu,
\hat\nu,\gamma,\hat\gamma$ in Theorem \ref{bettercenterbunch} are
constant, and choose $\s$ so that $\nu/\gamma<\s<
\min(1,\hat\gamma)$. We define the {\de center-unstable juliennes}
as \begin{equation} J^{cu}_n(x)=\bigcup_{y\in
W^c_{\s^n}(x)}f^{-n}(W^u_{\nu^n}(f^n(y)))=\bigcup_{y\in
W^c_{\s^n}(x)} J^u_n(y)
\end{equation}
Note that these are not solid juliennes, but laminae. However,
stable holonomy takes these laminae $J_n^{cu}(x)$ into
$J_{n-k}^{cu}(x')$ for all $n$, and some fixed $k$, so we have
attained a bounded distortion. If we saturate these sets by stable
leaves of length $\s^n$, we obtain a local base made of {\de
juliennes} $J_n(x)$, whose density points are preserved under stable
holonomy. This fulfills condition (2).
\par We must check now condition (1). As it was mentioned, it is
enough to see that the cubic system formed by the $\s^n$ $s-$ and
$u-$saturation of a center arc of length $\s^n$ is Vitali equivalent
to $J_n$ over $u$-saturated sets. We will see that both systems are
Vitali equivalent to a third one over essentially $u$-saturated
sets. Indeed, due to Lemma \ref{lema.abs.cont}, the cubic system is
easily seen to be equivalent to:
$$\bigcup_{y\in W^{cs}_{\s^n}(x)}J^u_n(y)$$
the rest of the proof consists in seeing that these new juliennes
constructed by first $\s^n$ $s$-saturating $W^c_{\s^n}(x)$, and then
saturating by $J^u_n(y)$, are Vitali equivalent to the ones obtained
above. This is not too difficult and can be found in full detail in
\cite{rhrhur1}, Proposition B. 9. With this, we have proved
Proposition \ref{jn vitali equivalent lebesgue}. With an analogous
procedure, we can obtain
\begin{proposition}
The set of Lebesgue density points of an essentially $s$- and
essentially $u$-saturated set is $u$-saturated.
\end{proposition}
With this, we get that the set of Lebesgue density points of an
essentially $s$- and essentially $u$-saturated set is $s$- and
$u$-saturated. The procedure consists in building another family of
juliennes, which are also equivalent to Lebesgue, but on essentially
$s$-saturated, sets. This can be attained by taking $\hat \s$ such
that $\hat\nu/\hat\gamma<\hat\s<\min(1,\gamma)$, and
$$J^{cs}_n(x)=\bigcup_{y\in
W^c_{\hat\s^n}(x)}f^{-n}(W^s_{\hat\nu^n}(f^n(y)))=\bigcup_{y\in
W^c_{\hat\s^n}(x)} J^s_n(y)$$ If we saturate by unstable leaves of
length $\hat\s^n$, we obtain a family of juliennes $J'_n$, which is
analogously seen to be equivalent to Lebesgue over essentially
$s$-saturated sets and $u$-saturated over essentially $u$-saturated
sets. 
\end{subsection}
\begin{subsection}{Some interesting corollaries}
\begin{proposition}
Let $f$ be a partially hyperbolic system. Given $\Sigma$ a
transversal to $E^s\oplus E^u$ there is a constant $C>0$ such that
if $A$ is an $su$-saturated set then
$$
\frac{1}{C}m(A)\leq m_{\Sigma}(A\cap\Sigma)\leq Cm(A)
$$
where $m_{\Sigma}$ is Lebesgue measure on $\Sigma$.
\end{proposition}

As a corollary of the above theorem we have:
\begin{corollary}\label{ergodictransversal}
Let $f$ be a bunched partially hyperbolic diffeomorphism. Let $P^c$
be a closed manifold everywhere tangent to $E^c$ such that
$f(P^c)=P^c$ and that every point can be joined to $P^c$ by an
$su$-path. If $f|P^c$ is ergodic then $f$ is ergodic. If $f|P^c$ is
stably ergodic then $f$ is stably ergodic.
\end{corollary}
Applying this corollary we get:
\begin{corollary}
Let $\phi$ be a Lebesgue measure preserving Anosov flow. Then either
$\phi_t$ is ergodic for every $t\in\R$ or $\phi$ is flow equivalent
to the suspension of an Anosov diffeomorphism by a constant function
$\omega$ and $t/\omega$ is rational.
\end{corollary}
Another corollary is:
\begin{corollary}
Let $f:M\to M$ be a stably ergodic diffeomorphism and let $g:N\to N$
be a volume preserving Anosov diffeomorphism, then, if $f\times g$
is bunched, using $M\times\{y\}$ as central foliations then it is
stably ergodic. Also, if $g:N\to N$ is partially hyperbolic with the
stable accessibility property, and $f\times g$ is bunched, using
$T_x M\times E^c_y$ as central space, then $f\times g$ is stably
ergodic.
\end{corollary}
What is interesting here is that one do not needs the accessibility
property for $f\times g$ to guaranty its stable ergodicity.
\begin{problem}
Are products of stably ergodic systems stably ergodic?
\end{problem}
\begin{problem}
Find an example of a volume preserving stably ergodic diffeomorphism
that is not robustly transitive, or either prove that there is no
such example.
\end{problem}

\end{subsection}

\begin{subsection}{The affine case}
Affine diffeomorphisms are always center bunched w.r.t. the
splitting $\g=\g^s\oplus\g^c\oplus\g^u$. Hence theorem
\ref{bettercenterbunch} always applies and hence $\overline{HB}=G$
implies ergodicity. Of course, for an affine diffeomorphism, one
does not need theorem \ref{bettercenterbunch} to prove ergodicity,
in fact it was already known by Dani, \cite{da1} \cite{da2}, that
$\overline{HB}=G$ implies the Kolmogorov property. But for
perturbations it is needed and hence, as a consequence of theorems
\ref{accessopenc1bund} and \ref{bettercenterbunch} the following is
also true
\begin{theorem}\cite{push3}
If an affine diffeomorphism has the accessibility property then it
is stably ergodic.
\end{theorem}
We say that an affine diffeomorphism is {\de stably ergodic among
left translations} if its perturbations by left translations are
also ergodic.
\begin{theorem}\cite{bresh}, \cite{push4}\label{stabergimpliesph}
If an affine diffeomorphism is stably ergodic among left
translations then it must be partially hyperbolic and
$\overline{HB}=G$.
\end{theorem}
Brezin and Shub have proven theorem \ref{stabergimpliesph} in the
semisimple and solvable cases and then Starkov proved it in full
generality. This motivated Pugh and Shub to formulate
\begin{problem}
Does $\overline{HB}=M$ implies stable ergodicity?, in other words,
is stable ergodicity among left translations enough for stable
ergodicity?
\end{problem}
So far, the only known examples of affine diffeomorphism being
stably ergodic but not having the accessibility property are the
automorphism of torus in \cite{frh}. On the other hand, Tahzibi has
asked
\begin{problem}
Does stable ergodicity for partially hyperbolic systems imply
essential accessibility when the central dimension is one? i.e. does
theorem \ref{stabergimpliesph} hold in this general context?
\end{problem}

\begin{problem}
Are the affine diffeomorphisms Bernoulli whenever $\overline{HB}=G$?
What happens with their perturbations?
\end{problem}
\end{subsection}
\begin{subsection}{Weak ergodicity}\label{weakergodic}
In the general case, the accessibility property gives a weak form of
ergodicity that in some cases it is used indeed to prove ergodicity.
We say that $f$ is {\de weakly ergodic} if almost every orbit is
dense. Obviously ergodicity implies weak ergodicity and in some
cases, weak ergodicity is enough to prove ergodicity.

Recall that $GC(x)=\bigcup_{n\in\Z}C(f^n(x))$.
\begin{proposition}\label{epsilonaccessibility}\cite{br1}, \cite{brpe}, \cite{budope}, \cite{dope}
Let $f$ be a partially hyperbolic system and assume that for
Lebesgue almost every point $x$, $GC(x)$ is $\epsilon$-dense. Then
Lebesgue almost every orbit is $\epsilon$-dense. In other words,
a.e. $\epsilon$-accessibility implies that a.e. orbit is
$\epsilon$-transitive.
\end{proposition}
We essentially take the proof from \cite{budope}.
\begin{proof}
Let $B$ be a ball of radius $\epsilon$. Let us say that a point $p$
is good if there is a neighborhood of $p$ such that a.e. point in
this neighborhood enters $B$. Then if we prove that a.e. point is
good we are done. Take $p$ such that $GC(p)$ is $\epsilon$-dense and
let us see that $p$ is good. By definition we have there is a
su-path $[z_0, \dots, z_k]$ with $z_0\in B$ and $z_k=f^N(p)$ for
some $N\in\Z$. It is clear that if $z_k$ is good then, as $f$ is a
diffeomorphism, $p$ will be good. We have obviously that $z_0$ is
good since it is itself in $B$. Let us see by induction that all the
$z_j$ are good. Assume $z_i$ is good and let us go to $z_{i+1}$. By
assumption $z_i$ has a neighborhood $N$ such that the orbit of a.e.
point in this neighborhood enters $B$. Take $S\subset N$ the whose
orbits enter $B$ and are forward and backward recurrent. By
Poincar\`e recurrence theorem we have that $S$ has full measure in
$N$. If $x\in S$ then the orbit of any point $y\in W^s(x)\cup
W^u(x)$ enters $B$. The absolute continuity of the foliations $W^s$
and $W^u$ means that the set
$$
\bigcup_{x\in S} W^s(x)\cup W^u(x)
$$
has full measure in the set
$$
\bigcup_{x\in N} W^s(x)\cup W^u(x)
$$
The latter is a neighborhood of $z_{i+1}$. Hence $z_{i+1}$ is good.
\end{proof}
As a corollary we get:
\begin{theorem}\label{accessimpliesweakerg}
If for Lebesgue a.e. $x$ $GC(x)$ is dense, $f$ is weakly ergodic.
\end{theorem}
Observe that the proof of proposition \ref{epsilonaccessibility}
only works for Lebesgue measure. Indeed, the main tools used are
Poincar\`e recurrence and absolute continuity of the strong
foliations. It would be interesting to have an analogous theorem for
other measures. See Subsection \ref{absolutecontinuity}.
\end{subsection}

\begin{subsection}{Ergodicity via Lyapunov exponents}\label{ergodiclyap}
A first corollary of theorem \ref{epsilonaccessibility} is that
accessibility plus local ergodicity somewhere implies ergodicity. In
\cite{budope} it is proven that negative central Lyapunov exponents
implies local ergodicity. But they also prove that all this
situation is stable under perturbations, in fact, they prove the
following:
\begin{theorem}\cite{budope}
Let $f$ be a volume preserving partially hyperbolic diffeomorphism
and assume that almost every accessibility class is dense. If the
exponents corresponding to the central bundle are all negative on a
positive measure set then $f$ is stably ergodic.
\end{theorem}
In this paper, Burns, Dolgopyat and Pesin pose the following
problem:
\begin{problem}
Is accessibility plus a.e. non-zero Lyapunov exponents enough for
ergodicity?
\end{problem}
In \cite{budope} the authors also used ideas from the papers
\cite{albovi}, \cite{bovi} where it is proven that with a dominated
splitting, see section \ref{rtransitivity} for a definition, and
some assumptions on the appearance of nonzero Lyapunov exponents one
obtains the existence of SRB measures. Also using this this type of
technics, Ali Tahzibi \cite{ta} found the following:
\begin{theorem}\cite{ta}
On $\T^n$, there is a stably ergodic diffeomorphism homotopic to an
Anosov diffeomorphism that admits a dominated splitting but has no
invariant hyperbolic subbundle.
\end{theorem}
Moreover he also proved the uniqueness of the SRB measure, that
exists by \cite{albovi}, for non-conservative perturbations. Hence
he proves in particular that partially hyperbolicity is not
necessary for proving stable ergodicity, see section
\ref{rtransitivity} for some related problems.
\end{subsection}
\end{section}

\begin{section}{Lyapunov exponents}\label{lyapunovexponent}
In this section we review some results relating Lyapunov exponents
and some other types of growth rates, with partially hyperbolic
systems.\par
Given a $C^1$ diffeomorphism $f:M\to M$, Oseledec's theorem asserts
the existence of some asymptotic directions with some asymptotic
growth rates corresponding to vectors in this directions for a.e.
point with respect to any invariant measure. Indeed, he proves the
following:
\begin{theorem}
Given an invariant measure $\mu$, for $\mu$-a.e. point $x$ there is
an invariant splitting $T_xM=E^1_x\oplus\dots\oplus E^{k(x)}_x$,
where $k(f(x))=k(x)$ and $Df_x(E^j_x)=E^j_{f(x)}$. There are also
invariant functions $\lambda_1(x)>\dots>\lambda_{k(x)}(x)$ such that
$$
\lim_{n\to\pm\infty}\frac{1}{n}\log|Df^n_x (v)|=\lambda_j(x)
$$
for $v\in E^j_x\setminus\{0\}$ and $j=1,
\dots, k(x)$. Moreover, if $1\leq j\neq l\leq k(x)$, then
$$
\lim_{n\to\pm\infty}\frac{1}{n}\log\angle
\left(E^j_{f^n(x)},E^l_{f^n(x)}\right)= 0.
$$
\end{theorem}
The numbers $\lambda_i$ are called the {\de Lyapunov exponents} and
the splitting is called the {\de Oseledec's splitting}. Notice that
if the measure is ergodic then the Lyapunov exponents and $k(x)$ are
constant a.e.\par
Observe that in the partially hyperbolic case, the Oseledec's
splitting refines the partially hyperbolic splitting. We call {\de
central Lyapunov exponents} (or central exponents) the Lyapunov
exponents corresponding to vectors in $E^c$, similarly for the
strong stable and strong unstable Lyapunov exponents. Observe that
the sum of the central Lyapunov exponents, counted with
multiplicities equals:
$$
\sum_j\lambda_j(x)\dim
E^j_x=\lim_{n\to\infty}\frac{1}{n}\sum_{k=0}^{n-1}\log Jac\left(
D_{f^k(x)}f|E^c\right)
$$
where the sum ranges over all $j$ such that $E^j_x\subset E^c_x$. We
will be interested in the integrated sum of central Lyapunov
exponents:
$$
\int\sum_j\lambda_j(x)\dim E^j_xd\mu
$$
which equals
$$
\int\log Jac\left( D_xf|E^c\right)d\mu
$$
and, as always,  in most cases the invariant measure is Lebesgue
measure. Notice that when the measure is ergodic the integrated sum
of central Lyapunov exponents equals a.e. the sum of central
Lyapunov exponents.

\begin{subsection}{Removing zero central exponents}
Let us begin with Shub--Wilkinson result about the approximation by
non-zero central Lyapunov exponents.
\begin{theorem}\cite{shwi2}
There is an open set of partially hyperbolic volume preserving
Bernoulli diffeomorphism with non-zero Lyapunov exponents on $\T^3$.
\end{theorem}
The idea is to get stably ergodicity and non-zero exponents, then
Pesin's results give the Bernoulli property \cite{pe}. They begin
with a linear map $f$ of $\T^3=\T^2\times\T$ associated with the
matrix
\begin{displaymath}
B=\left(\begin{array}{cc}
A           & 0\\
{\mathbf w_0} & 1
\end{array}\right),
\end{displaymath}
where where ${\mathbf w_0}\in\Z^2\setminus\{0\}$ and $A\in
GL(2,\Z)$, for example,
\begin{displaymath}
A=\left(\begin{array}{cc}
2 & 1\\
1 & 1
\end{array}\right).
\end{displaymath}
and ${\mathbf w_0}=(1,1)$. They embed $f$ into a two parameter
family $f_{a,b}$ such that $f_{0,0}=f$. With the parameter $b$ they
will guarantee accessibility and hence ergodicity and with the
parameter $a$ they will get non-zero exponents. It is important in
order to get non-zero exponents in their construction that ${\mathbf
w_0}\neq 0$. Although a posteriori one can arrage the example to get
one close to Anosov times identity.

Let us see how it works. Observe that for a partially hyperbolic
diffeomorphism with central dimension $1$ the sum of the central
exponents is just the central exponent. They look for a family
$f_{a,b}=g_a\circ h_b$ where $h_b$ is a skew product over $A$ having
the accessibility property and $g_a$ is a perturbation along the
unstable direction that will make the strong unstable exponent
decrease without touching the strong stable exponent. Let us denote
$\x=(x,y)$, then
$$
h_b(\x,z)=\left(A\x,z+{\mathbf w_0}\cdot\x+b\varphi(\x)\right)
$$
where $\varphi:\T^2\to\R$ is suitably chosen, for example
$\varphi(x,y)=\sin(2\pi y)$. We shall see later what is meant by
suitable in general. And
$$
g_a(\x,z)=(\x+a\psi(z){\mathbf v_0},z)
$$
where $\psi:\T\to\R$ is also suitably chosen, for example
$\psi(z)=\sin(2\pi z)$, and ${\mathbf v_0}$ is a unit vector in the
unstable direction, for instance ${\mathbf
v_0}=\left((1+\sqrt{5})/2,1\right)$. Observe that $h_b$ and $g_a$
are volume preserving for every $a,b$.\par
Let us see what is suitable for $\varphi$. The idea here is to apply
the following criterion to guarantee accessibility for skew
products, see \cite{brpe, buwi1}. Let $\theta:\T^2\to\T$ be a smooth
map and define the skew product
$h_{\theta}:\T^2\times\T\to\T^2\times\T$ by
$h(\x,z)=(A\x,z+\theta(\x))$ where $A\in GL(2,\Z)$ is a hyperbolic
matrix. Then $h$ is partially hyperbolic and has the accessibility
property if and only if there is no solution to the cohomological
equation
\begin{eqnarray}\label{cocycleeq}
k\theta(\x)=\Phi(A\x)-\Phi(\x)+c
\end{eqnarray}
for $\Phi\in C^0(\T^2,\T)$, $c\in\T$, where $k=\det(A-Id)$.

Observe that $\theta$ may be written as $\theta(\x)={\mathbf
w_0}\cdot\x+\varphi(\x)$ where ${\mathbf w_0}\in\Z^2$ and $\varphi$
is homotopic to constant and hence may be seen as a function
$\varphi:\T^2\to\R$. Similarly, the unknown $\Phi$ may be written as
$\Phi(\x)={\mathbf w_1}\cdot\x+k\eta(\x)$ where ${\mathbf
w_1}\in\Z^2$ and $\eta$ is homotopic to constant and hence also
$\eta:\T^2\to\R$. Observe that in order that the cohomological
equation \ref{cocycleeq} have a solution, it is necessary that
$(A^t-Id){\mathbf w_1}=k{\mathbf w_0}$. And this has solution
${\mathbf w_1}\in\Z^2$ since $k$ is exactly
$\det(A^t-Id)=\det(A-Id)$. Hence taking ${\mathbf w_1}$ this
solution, the cohomological equation \ref{cocycleeq} transforms into
\begin{eqnarray}\label{cocycleeq2}
\varphi(\x)=\eta(A\x)-\eta(\x)+c
\end{eqnarray}
for $\eta\in C^0(\T^2,\R)$ and $c\in\R$. Finally, notice that for
the cohomological equation \ref{cocycleeq2} to have solution it is
necessary that the average of $\varphi$ over any invariant measure
be $c$. So that any function $\varphi$ whose average at two
different periodic points differs will be good to guarantee
accessibility and hence $\varphi$ will be suitable.\par
So, for any $b\neq 0$, $f_{0,b}=h_b$ has the accessibility property
and hence it belongs to an open set of ergodic diffeomorphisms. So,
for any fixed $b$ one can move the parameter $a$ a little still
having ergodicity. Let us see what happens when we move $a$,let us
first look at the invariant foliations for $a=0$. It follows that
the central, the center-stable and center-unstable foliations are
still the same as for $b=0$. For the stable foliation we have that
each leaf is the graph of a function from the stable manifold of $A$
to the circle and it is invariant under translations of the form
$(\x,z)\to (\x,z+z_0)$ for every $z_0\in\T$. Now, for any $a$, $g_a$
preserves the center-unstable foliation by planes. And hence
$f_{a,b}$ also preserves the center-unstable foliation by
planes.\par
Let us see that the strong stable Lyapunov exponent for $f_{a,b}$ is
the same as for $f_{0,0}$. To that end one can work in the universal
covering and hence there is a well defined linear projection
$\pi^s:\R^3\to E^s_B$ along the center-unstable planes into the
stable manifold of the linear one such that $\pi^s\circ
f=B\circ\pi^s$. This is because $f$ preserves the center-unstable
foliation by planes. Hence $\pi^s(D_xf(v))=\sigma^{-1}\pi^s(v)$ for
any vector $v$, where $\sigma$ is the unstable eigenvalue. So, if we
put the sup-norm $|v|=\max\{|v^s|,|v^{cu}|\}$ where $v=v^s+v^{cu}$
w.r.t. the splitting for the linear, we get that for any vector $V$
in the stable direction of $f_{a,b}$, $|V|=|V^s|$ since the
invariant spaces for the perturbation are close to the linear ones.
Hence it follows that
$$
\left|D_xf^{-n}(V)\right|=\left|\pi^s\left(D_xf^{-n}(V)\right)\right|=\sigma^n|V^s|=\sigma^n|V|
$$
and hence
$$
\frac{1}{n}\log\left|D_xf^{-n}|E^s_{a,b}\right|=\log\sigma
$$
for any $n$ which implies that the strong stable Lyapunov exponents
coincide. Now, as $f_{a,b}$ is volume preserving it follows that the
sum of the strong unstable exponent and the center exponent should
equal $\log\sigma$. Hence, if the strong unstable exponent decrease,
 the central exponent should be positive and that is what should be
computed.\par
The idea again is to use that $f_{a,b}$ is a perturbation so that we
may write a vector-field $V^u(\x,z)=({\mathbf v_0},u_{a,b}(\x,z))$
where ${\mathbf v_0}$ is the vector in the unstable direction, and
$u_{a,b}$ is uniquely defined by this requirement. Observe that
$u_{0,0}=\frac{{\mathbf v_0}\cdot{\mathbf w_0}}{\sigma-1}\neq 0$.
Then after some computation they obtain that the strong unstable
exponent is
$$
\log\sigma-\int_{\T^3}\log[1-a\psi'(w)u_{a,b}(w)]dw
$$
Now, they see how this number varies and they get that it decreases.
That is made by some heavy computation calculating its derivatives.

The idea is somehow that when $a$ is non zero, the unstable
distribution will have some component on the center bundle and hence
it will force the dynamics on the strong unstable to slow down.\par
In the $C^1$ setting, Alexandre Baraviera and Christian Bonatti were
able to push this technic to a more general context proving:
\begin{theorem}\cite{babo}\label{barabona}
For a $C^1$ open and dense set of volume preserving diffeomorphism
$f:M\to M$ admits a dominated splitting $TM=E_1\oplus\dots E_k$,
$k>1$ such that the integrated sum of Lyapunov exponents on  $E_i$
are non zero for $i=1, \dots, k$, i.e. $\int_M\log Jac
\left(D_xf|E_i\right)dx$.
\end{theorem}
For the definition of dominated splitting see section
\ref{rtransitivity}.

\end{subsection}

\begin{subsection}{$C^1$ genericity and Lyapunov exponents.}
In \cite{mane3} it appeared a sketch of a proof of the following:
\begin{theorem}\cite{boc, mane3}
For a $C^1$-generic area preserving diffeomorphism of surfaces $f$,
either $f$ has a.e. zero exponents or $f$ is Anosov
\end{theorem}
In particular if the surface is not a torus then the generic $f$ has
zero exponent. In \cite{boc}, Jairo Bochi fixed and completed the
proof in \cite{mane3} and then in \cite{bocvi}, jointly with Marcelo
Viana, they generalized the result to any dimension but loosing a
little of strength. In section \ref{rtransitivity} the reader may
find a definition of dominated splitting. We say that the Oseledecs
splitting is {\de dominated at $x$} if it extend to a dominated
splitting on the closure of the orbit of $x$.
\begin{theorem}\label{bochiviana}
For a $C^1$-generic volume preserving diffeomorphism $f$ of a
manifold $M$, for a.e. $x\in M$ the Oseledecs splitting of $f$ is
either trivial or dominated at $x$.
\end{theorem}
The techniques here involve the use of some tower constructions and
some perturbations that resemble the Shub-Wilkinson case in order to
slow down an exponent when there is no domination, but here of
course the slowdown will be not continuous. In fact, they deduce
theorem \ref{bochiviana} from the following
\begin{theorem}
Let $f_0$ be a $C^1$ volume preserving diffeomorphism such that the
map
$$
f\to (LE_1(f), . . . , LE_{d-1}(f))
$$
is continuous at $f=f_0$. Then for almost every $x\in M$, the
Oseledecs splitting of $f$ is either dominated or trivial at x.
\end{theorem}
Here $M$ is $d$-dimensional and $LE_ i(f)$ is the integrated sum of
the first $i$ Lyapunov exponents.

Of course it will be interesting to get a dominated slitting on the
whole manifold. Moreover, maybe the following is true:
\begin{problem}
$C^1$ generically among volume preserving diffeomorphisms, either
$f$ has a.e zero Lyapunov exponents or $f$ is ergodic.
\end{problem}
They also have some counterpart for the symplectic case stating that
a $C^1$ generic symplectic diffeomorphism is either Anosov or has at
least one (necessarily double) zero Lyapunov exponent. But a
complete counterpart is still open, that is, they ask if one can get
in fact partial hyperbolicity along the orbit of $x$. They put the
following problem:
\begin{problem}
Is it true that the Oseledecs splitting of generic symplectic $C^1$
diffeomorphisms is either trivial or partially hyperbolic at almost
every point?
\end{problem}
Indeed one can ask to get either zero exponents or partial
hyperbolicity. On the other hand, in the partially hyperbolic
setting, one can put together theorems \ref{bochiviana},
\ref{barabona}, \ref{dolgowilkinson} and \ref{accessimpliesweakerg}
plus some trick to bypass some absolute continuity and get that a
$C^1$ generic volume preserving, partially hyperbolic system has a
globally defined dominated splitting that coincides a.e. with the
Oseledecs splitting and also has non-zero integrated Lyapunov
exponents. So, on each bundle, the dynamics is asymptotically
conformal for a.e. point.
\end{subsection}

\begin{subsection}{Dynamical growth rates}
The Lyapunov exponents, as we have seen, measure the growth rate of
the derivative along some directions. There are others types of
asymptotic growth rates that can be defined, and they are all
typically related in some way. Let $f:M\to M$ be a diffeomorphism
and assume it leaves a foliation $\F$ invariant, a priori, as
always, with smooth leaves of dimension $d$ and tangent to a
continuous sub-bundle $E\subset TM$. We shall also assume  that $f$
expands $\F$. Let us define the {\de dynamical volume growth} of $f$
on $\F$ as
$$
vg_r(x)=\liminf_{n\to\infty}\frac{1}{n}\log\vol\left(f^n\left(\F_r(x)\right))\right)
$$
Where $\vol(A)$ is the volume of $A$ on the corresponding leaf
volume. This dynamical growth of volume was treated in the general
case by Newhouse \cite{ne2} and Yomdin \cite{yo}, in this last paper
is also proven the entropy conjecture for $C^{\infty}$ maps, that
is, the hyperbolicity of the action of $f$ in homology is a lower
bound for the entropy of $f$.

Then we define the {\de dynamical homological growth} of $f$ on $\F$
as the current defined in the following way, \cite{rusu}, given a
$d-$form $\omega$:
$$
C_{r,x}(\omega)=\lim_{n\to\infty}\frac{1}{\vol\left(f^n\left(\F_r(x)\right))\right)}\int_{\F_r(x)}f_*^n\omega
$$
Recall that a {\de current} is an element of the dual of the
differential forms. It turns out that $C$ is a closed current (it
vanishes on exact forms) and hence it defines an homology class. In
\cite{push3} Pugh and Shub asked:
\begin{problem}\label{curr}
Do the strong stable and unstable manifolds represent non-trivial
homology classes in the homology of M?.
\end{problem}
In the Anosov case the answer is positive, see \cite{ruwi}. In
\cite{ru}, Ruelle also defines some transverse measure to $\F$,
$\rho$, associated to some cocycles and such that
$f_*\rho=\lambda\rho$ for some $\lambda>1$. In particular when the
cocycle is trivial he gets a transverse invariant measure $\rho_0$
with $f_*\rho_0=\lambda_0\rho_0$.\par
All this quantities typically exists for almost every point and do
not depend on $r$. They are all related when the foliation is
absolutely continuous. This is the subject of a work of Saghin and
Xia that uses this relations to build as a corollary some
non-absolutely continuous invariant foliations.

\end{subsection}

\begin{subsection}{Lyapunov exponents and uniform hyperbolicity}
It is known that if a continuous map $f:M\to M$ is uniquely ergodic
then the Birkhoff averages converge uniformly. With essentially the
same proof one can prove that for a function $\varphi$, if the
integral of $\varphi$ w.r.t. any ergodic invariant measure is less
than a constant $C$, then the Birkhoff averages should be less than
$C$, i.e. $1/n\sum_{k=0}^{n-1}\varphi(f^k(x))<C$ for every $x$ and
$n\geq N_0$. In particular, if the integral w.r.t. any ergodic
invariant measure is always the same, then the Birkhoff averages
converge uniformly. The results in this subsection can be found in
\cite{frh2}. The reader may found results related to the ones in
this subsection in \cite{alarsa, ca1, ca2, caluri, sc}, we would
like to thank Yongluo Cao for putting these references into our
attention.

The next proposition says that for the multiplicative case we have
essentially the same phenomenon.

\begin{proposition}\label{subaditiv}
Let $f:X\to X$ be a continuous map of a compact metric space. Let
$a_n:X\to\R$, $n\geq 0$ be a sequence of continuous functions such
that $a_{n+k}(x)\leq a_n(f^k(x))+a_k(x)$ for every $x\in X$,
$n,k\geq 0$ and such that there is a sequence of continuous
functions $b_n$, $n\geq 0$ satisfying $a_n(x)\leq
a_n(f^k(x))+a_k(x)+b_k(f^n(x))$ for every $x\in X$, $n,k\geq 0$. If
$$
\inf_n \frac{1}{n}\int_X a_nd\mu<0
$$
for every ergodic $f$-invariant measure, then there is $N\geq 0$
such that $a_N(x)<0$ for every $x\in X$.
\end{proposition}

An interesting case where we will apply the proposition is the case
when $a_n(x)=\log |D_xf^n|E|$ and $b_n(x)=\log m(D_xf^n|E)$ getting
the following corollary. A regular $C^1$ map is a map whose
derivative is invertible at each point.
\begin{corollary}
Let $f:M\to M$ be a regular $C^1$ map and $\Lambda$ a compact
invariant set. Assume $f$ leaves a continuous bundle $E$ over
$\Lambda$ invariant. If the Lyapunov exponents of the restriction of
$Df$ to $E$ are all negative (positive) for every ergodic invariant
measure, then $Df$ contracts (expands) $E$ uniformly.
\end{corollary}
Also, using the fact that hyperbolic measures (measures with nonzero
Lyapunov exponents) are sent to hyperbolic measures by H\"older
continuous conjugacies, we have the following:
\begin{corollary}
Let $f:M\to M$ be a diffeomorphism and $g:N\to N$ be a
$C^{1+\mbox{\tiny{H\"older}}}$ diffeomorphism. Let $\Lambda$ be a
transitive hyperbolic set for $f$ and assume there is a H\"older
continuous homeomorphism $h:U\to V$ from a neighborhood $U$ of
$\Lambda$ onto $V\subset N$ such that $h\circ f=g\circ h$. Let us
assume that $g$ leaves a continuous splitting $TM=E_1\oplus E_2$
over $h(\Lambda)=\Lambda_g$ invariant, and that it coincides with
the Lyapunov (stable$\oplus$unstable) splitting for some
(necessarily hyperbolic) $g$-invariant measure. Then $\Lambda_g$ is
a hyperbolic set for $g$.
\end{corollary}
Also related to regularity of the invariant distributions we have
the following:
\begin{corollary}
Let $g$ be a $C^k$ Anosov diffeomorphism, and assume it preserves a
continuous splitting $TM=E_1\oplus E_2$ (not necessarily the
hyperbolic splitting). Given a periodic point $p$, let us call
$\chi_1^+(p)$ the biggest Lyapunov exponent of the restriction of
$Df$ to $E_1$, $\chi_2^+(p)$ the biggest Lyapunov exponent of the
restriction of $Df$ to $E_2$ and $\chi_2^-(p)$ the smallest Lyapunov
exponent of the restriction of $Df$ to $E_2$. If there is a constant
$c<0$ such that $\chi_1^+(p)-\chi_2^-(p)<c<0$ and
$\chi_1^+(p)+r\chi_2^+(p)-\chi_2^-(p)<c<0$, where $r\geq 1$, for
every periodic point $p$ then there is a $C^s$ foliation tangent to
$E_1$ where $s=\min\{k-1,r\}$.
\end{corollary}

\end{subsection}

\end{section}

\begin{section}{Integrability of the central
distribution}\label{integrability}

The integrability of the central distribution is one of the more
striking problems in the study of partial hyperbolic systems.
Indeed, there are quite few new results towards such integrability.
In general, given a plane field $E\subset TM$, there are two
possible obstructions to the integrability of $E$.
\begin{enumerate}[i)]
\item One obstruction  is that $E$ does not satisfy the  Froebenius bracket condition.
\item The other one is the lack of differentiability of the bundle itself.
\end{enumerate}
In the partially hyperbolic setting, let us mention that although
there are examples of non integrable central distributions, in this
examples the problem is Froebenius bracket condition and not the
differentiability, see subsubsection \ref{examplenoninteg}.
Moreover, subsection \ref{smoothcase} suggest that the Froebenius
part of the problem is intimately related to bunching, that is, if
$f$ satisfies some bunching then $E^{cs}$ and $E^{cu}$ should be
``involutive".

\begin{subsection}{The smooth case}\label{smoothcase}
In this first section we shall deal with the first reason of non
integrability mentioned above. First let us see a positive result
and then an example of non-integrability.

In \cite{buwi2} it is proven the integrability of the center
distribution when some bunching condition is available.
\begin{theorem}\cite{buwi2}
If $E^{cu}$ is smooth and the bunching condition
$\hat{\nu}<\hat{\gamma}^2$ holds then $E^{cu}$ is integrable,
analogously, if $E^{cs}$ is smooth then it is integrable whenever
$\nu<\gamma^2$. Finally, if both $E^{cs}$ and $E^{cu}$ are smooth
and both bunching conditions hold, $E^c$ is integrable.
\end{theorem}
In \cite{buwi2} the authors give a geometric proof. Their proof
essentially follows the lines of the proof of Froebenius' theorem.
We shall present two proofs of the theorem. First a wrong proof that
uses Froebenius' theorem itself. Second, a correction of the first
proof.

\begin{proof}{\bf First (wrong) proof.}
Let as see the $E^{cs}$ case, the other case is analogous. By
Froebenius' theorem, it has to be proven that whenever $X$ and $Y$
are two vector-fields tangent to $E^{cs}$ their Lie bracket $[X,Y]$
is also tangent to $E^{cs}$. If we denote the norm by $|\cdot|$, we
have that
$$
|D_pf^n\left([X,Y]\right)|=|\left[D_pf^n(X),D_pf^n(Y)\right]|\leq
C|D_pf^n(X)||D_pf^n(Y)|\leq C\hat{\gamma}_n(p)^{-2}|X||Y|.
$$
since $X,Y\in E^{cs}$. On the other hand, recall that $E^{cs}$ may
be characterized as the vectors where $|D_pf^n(Y)|\hat{\nu}_n(p)\to
0$. Finally, as $\hat{\nu}<\hat{\gamma}^2$ we get that
$C\hat{\gamma}_n(p)^{-2}\hat{\nu}_n(p)|X||Y|\to 0$ and we are done.

What is wrong in this proof is that we are tacitly using that
$|[X,Y]|\leq C|X||Y|$ for some constant $C$. That is false, what you
have is that $|[X,Y]|\leq C|X|_{C^1}|Y|_{C^1}$, but this inequality
is much more complicated to deal with.
\end{proof}
\begin{proof}{\bf Second (corrected) proof.}
We use again Froebenius' theorem, let $X$ and $Y$ be two
vector-fields tangent to $E^{cs}$ and let us see that  their Lie
bracket $[X,Y]$ is in $E^{cs}$. Given $p\in M$ take $z\in\omega(p)$,
and $\eta^1,\dots,\eta^u$, $u$ linearly independent 1-forms defining
$E^s\oplus E^c$ in a neighborhood of $z$, that is, $E^s\oplus E^c$
is the intersection of the kernels of $\eta^i$. Take $n_j\to\infty$
such that $p_j=f^{n_j}(p)\to z$ and call
$$
v_j=\frac{D_pf^{n_j}[X,Y](p)}{\left|D_pf^{n_j}[X,Y](p)\right|}
$$
and we may assume that $v_j\to v$. Now, if $[X,Y](p)\notin
E^s_p\oplus E^c_p$, then, $v_j$ looses its center-stable component
and its unstable component persists, so that $v\in E^u_z$ and $v\neq
0$. Thus we get that there is $i$ such that $\eta^i_z(v)\neq 0$ and
thus $|\eta^i_{p_j}(v_j)|=|\eta^i_{p_j}(v_j^u)|>c>0$, where
$v_j=v_j^{cs}+v_j^u$. Let us call $X_{n_j}(x)=D_xf^{n_j}X(x)$ and
$Y_{n_j}(x)=D_xf^{n_j}Y(x)$. We may assume that $p_j$ are in the
neighborhood of $z$ so that, on one hand we have

\begin{eqnarray*}
|\eta^i_{p_j}(D_pf^{n_j}[X,Y](p))|&=&
|\partial_{X_{n_j}}\eta^i_{p_j}(Y_{n_j})-\partial_{Y_{n_j}}\eta^i_{p_j}(X_{n_j})-d\eta^i_{p_j}(X_{n_j}(p),Y_{n_j}(p))|
\\&=&
|d\eta^i_{p_j}(D_pf^{n_j}X(p),D_pf^{n_j}Y(p))|\\
&\leq&C|D_pf^{n_j}X(p)| |D_pf^{n_j}Y(p)|\\&\leq&
C\hat{\gamma}_{n_j}(p)^{-2}|X||Y|
\end{eqnarray*}
and on the other hand we have that
\begin{eqnarray*}
|\eta^i_{p_j}(D_pf^{n_j}[X,Y](p))|&=&|D_pf^{n_j}[X,Y](p)||\eta^i_{p_j}(v_j)|\\
&\geq&|D_pf^{n_j}[X,Y]^u(p)||\eta^i_{p_j}(v_j^u)|\geq
c\hat{\nu}_{n_j}(p)^{-1}|[X,Y]^u|
\end{eqnarray*}
where $[X,Y]=[X,Y]^{cs}+[X,Y]^u$. Thus, since $[X,Y]^u\neq 0$, we
get that $\hat{\nu}_{n_j}(p)^{-1}\leq C\hat{\gamma}_{n_j}(p)^{-2}$
which contradicts the bunching condition.

\end{proof}

It would be interesting to prove the integrability of the central
distribution but only assuming the differentiability of $E^c$ and
bunching. In fact, this result lead Keith Burns and Amie Wilkinson
to the formulation of the following problem:
\begin{problem}\label{burnswilkinteg}
If $f$ is center-bunched then the center bundle is integrable.
\end{problem}
Related to this problem we want to mention the weak-integrability
notion defined in \cite{brbuiv}. A bundle $E\subset TM$ is said to
be weakly integrable if for every point $x\in M$ there is a complete
manifold $x\in W(x)$ tangent to $E$. In \cite{brbuiv} it is proven
that if the center bundle is one dimensional then $E^{cs}$, $E^{cu}$
and $E^c$ are weakly integrable. Also they prove that if two
partially hyperbolic diffeomorphisms are homotopic through a path of
partially hyperbolic diffeomorphisms and one of them has $E^{cs}$
weakly integrable then so does the other. Let us put a problem
weaker than problem \ref{burnswilkinteg}:
\begin{problem}
If $f$ is center-bunched, is the center bundle weakly-integrable?.
\end{problem}
Observe that if the bundle $E$ is smooth then weak-integrability and
integrability coincide.

\begin{subsubsection}{An example of non-integrability.}\label{examplenoninteg}
Let $\He$ be the Heisenberg group as in subsection \ref{sacksteder}.
Then $f$ will be essentially A. Borel's construction of an Anosov
diffeomorphism on a quotient of $\He^2=\He\times\He$ that appeared
in Smale's paper \cite{sm}. The Lie algebra of $\He^2$ is
$\he\oplus\he$, where the bracket is defined componentwise. Let
$\lambda=2+\sqrt{3}$, then $\lambda^{-1}=2-\sqrt{3}$. Given $\alpha,
\beta\in\Z$ let us define the automorphism of $\he\oplus\he$
$$
f(x_1,x_2,y,a_1,a_2,b)=(\lambda^{\alpha}x_1, \lambda^{\beta}x_2,
\lambda^{\alpha+\beta}y, \lambda^{-\alpha}a_1, \lambda^{-\beta}a_2,
\lambda^{-\alpha-\beta}b).
$$
It is well defined and hence it defines an automorphism $F$ of
$\He^2$. It is partially hyperbolic if $\alpha$ or $\beta$ are not
zero and hyperbolic if $\alpha, \beta, \alpha+\beta$ are non-zero.
Now the whole theme is to find a lattice $\Gamma$ of $\He^2$ that is
$F-$invariant and cocompact, i.e with $\He^2/\Gamma$ compact. But
before doing this let us identify the invariant bundles, and the
invariant foliations already in $\he\oplus\he$, recall that
$\he_1=[\he,\he]$ and that $[\he_1,\he]=0$. Let us assume that
$\alpha+\beta>\beta\geq\alpha$. Then we can take
$E^u=\he_1\oplus\{0\}$ to be the strong unstable direction and
$E^s=\{0\}\oplus\he_1$ to be the strong stable one. Then take $E^c$
to be the space spanned by $X_1=(1,0,0,0,0,0), X_2=(0,1,0,0,0,0),
A_1=(0,0,0,1,0,0)$ and $A_2=(0,0,0,0,1,0)$, that is, the orthogonal
complement to $\he_1\oplus\he_1=E^u\oplus E^s$. Then observe that
$E^c$ is not integrable since $[E^c,E^c]=\he\oplus\he$ and hence it
is not involutive. In fact $0\neq Y=[X_1,X_2]\in\he_1\oplus\{0\}$
and $0\neq B=[A_1,A_2]\in\{0\}\oplus\he_1$. Observe that also we
could also take $E^s=\{0\}\oplus\he$ and then $E^c$ the orthogonal
complement of $\he_1\oplus\he$ and then $E^c\subset\he\oplus\{0\}$
and $[E^c,E^c]=\he\oplus\{0\}$ and it is still not integrable.

Let us see now how to build an $F-$invariant, cocompact lattice. We
take this part from \cite{ausch}. Take, $a\in\Z$ such that $\lambda$
be a root of $x^2+2ax+1$ and define the following basis for
$\He\times\He$,
$$
\left\{X_1+A_1, \sqrt{a^2-1}(X_1-A_1), X_2+A_2,
\sqrt{a^2-1}(X_2-A_2), Y+B, \sqrt{a^2-1}(Y-B)\right\}
$$
and call the vectors $E_1,\dots, E_6$ respectively, then it follows
that $[E_1,E_3]=E_5, [E_1,E_4]=E_6, [E_2,E_3]=E_6$ and
$[E_2,E_4]=(a^2-1)E_5$ and all other possible combinations are
either zero or the opposite of the existing ones. So taking the
integer lattice $L\subset\he\oplus\he$ formed by the integer
combination of $E_1,\dots E_6$, and $\Gamma=\exp(L)$ we get that
$\Gamma$ is a discrete subgroup of $\He^2$ and that $\He^2/\Gamma$
is a compact nilmanifold. Let us see that $F$ defines an
automorphism of $\He/\Gamma$, for this we need to see that $f(L)=L$
or, what is the same, that the matrix associated to $f$ in the basis
$\{E_i\}$ is an integer matrix. Then, if we take the two by two
matrix
\begin{displaymath}
C=\left(\begin{array}{cc}
-a & a^2-1 \\
 1 & -a
\end{array}\right)
\end{displaymath}
then the assosiated matrix to $f$ is:
\begin{displaymath}
\left(\begin{array}{ccc}
 C^{\alpha} & 0         & 0\\
 0          & C^{\beta} & 0\\
 0          & 0         & C^{\alpha+\beta}
\end{array}\right)
\end{displaymath}
In the case $\lambda=2+\sqrt{3}$, $a=-2$ and if we take $\alpha=1$
and $\beta=2$ we obtain
\begin{displaymath}
\left(\begin{array}{cccccc}
 2 & 3 & 0  & 0  & 0  & 0\\
 1 & 2 & 0  & 0  & 0  & 0\\
 0 & 0 & 7  & 12 & 0  & 0\\
 0 & 0 & 4  & 7  & 0  & 0\\
 0 & 0 & 0  & 0  & 26 & 45\\
 0 & 0 & 0  & 0  & 15 & 26
\end{array}\right)
\end{displaymath}

Observe that using some variation of theorem \ref{accessopenc1bund}
it can be proven that for any perturbation of $F$, the center-stable
and the center-unstable bundles are not integrable. In fact, one can
prove that any two points can be joined by central curves. Also it
follows the stability of the non integrability using the notion of
weak-integrability in \cite{brbuiv}. In fact, neither $E^{cs}$, nor
$E^{cu}$ are weakly-integrable and this is a widely stable property.

\end{subsubsection}
\end{subsection}

\begin{subsection}{Some special cases}
In \cite{rhrhur2} we prove the unique integrability of the central
distribution when the diffeomorphism $f$ is transitive, the manifold
is of dimension three  and there are no periodic points at all. When
the center-stable or the center-unstable is integrable, if the
fundamental group of $M$ is abelian, the isomorphism induced by $f$
in the homology is partially hyperbolic, see \cite{brbuiv} and
subsection \ref{growthcurves}. This implies that, in particular, if
$M=\S^3$ then $Per(f)\neq\emptyset$. However, it is not known
whether  $\S^3$ supports a partially hyperbolic diffeomorphism.

\begin{problem}
Does $\S^3$ support a partially hyperbolic diffeomorphism?
\end{problem}

We think the answer is negative. For such a diffeomorphism, neither
$E^{cu}$ nor $E^{cs}$ could be integrable (see
\cite{dipuur,brbuiv}).

\begin{subsubsection}{Quasi-isometric foliations and integrability}
A foliation $\W$ of a simply connected riemannian manifold is said
to be quasi-isometric if there are constants $a$ and $b$ such that
whenever $x$ and $y$ are in the same leaf $d_W(x,y)\leq ad(x,y)+b$,
where $d_W(x,y)$ denotes the distance on the leaf. In \cite{br3}
Brin prove the following:
\begin{theorem}
Let $f$ be a partially hyperbolic diffeomorphism of a compact
riemannian manifold $M$. Suppose the unstable foliation $\W^u$ of
$f$ is quasi-isometric in the universal cover $\tilde{M}$. Then the
distribution $E^{cs}$ is locally uniquely integrable.
\end{theorem}
Of course if $\W^s$ is quasi-isometric in $\tilde{M}$ then $E^{cu}$
is locally uniquely integrable, and if both are quasi-isometric in
$\tilde{M}$ then $E^c$ is locally uniquely integrable. Although the
unstable foliation do not need to be quasi-isometric in the
universal covering, for example for the geodesic flow on a
hyperbolic surface, it is quasi-isometric in some interesting
examples.
\begin{proposition}\cite{br3}
Let $\W$ be a $k-$dimensional foliation on $\T^m$. Suppose there is
a codimension $k$ plane $A$ such that $TW\cap A=\{0\}$. Then the the
lift of $\W$ is quasi-isometric in the universal cover $\R^m$.
\end{proposition}
It would be interesting to weaken somehow the hypothesis of the
proposition, at least to a more topological type of hypothesis. In
\cite{br3} Brin ask the following:
\begin{problem}
If $f$ is a partially hyperbolic diffeomorphism of $\T^3$, are the
stable and unstable foliations necessarily quasi-isometric?
\end{problem}
Of course we may ask the same for $\T^n$.
\end{subsubsection}

For other cases of unique integrability of the center bundle see
subsection \ref{transbowi}.

\end{subsection}

\begin{subsection}{Plaque expansive}
Given a partially hyperbolic diffeomorphism $f:M\to M$ having a
center foliation $\F^c$, we define a $\delta-$pseudo-orbit
respecting the central plaques to be a sequence $x_n$, $n\in\Z$,
such that $f(x_n)\in\F^c_{\epsilon}(x_{n+1})$. We say that $f$ is
plaque expansive at $\F^c$ if there is an $\epsilon>0$ such that if
$x_n$ and $y_n$ are $\epsilon-$pseudo-orbits preserving the central
plaques and $d(x_n,y_n)<\epsilon$ for every $n\in\Z$ then
$x_0\in\F^c_{\epsilon}(y_0)$. The main reference for plaque
expansivity is still \cite{hipush}.

If $\F^c$ is a $C^1$ foliation it is plaque expansive, \cite{hipush}
otherwise it is only known in some cases.
\begin{problem}
Are the central foliations always plaque expansive? What about when
the strong foliations are quasi-isometric?
\end{problem}
One of the main consequences of plaque expansiveness is the
following:
\begin{theorem}\cite{hipush}
Let $f:M\to M$ be a plaque expansive partially hyperbolic
diffeomorphism then there is a neighborhood of $f$, $U$, such that
if $g\in U$ then $g$ leave invariant a plaque expansive center
foliation $\F^c_g$ and there is an homeomorphism $h:M\to M$ such
that $h\left(\F_f(x)\right)=\F_g\left(h(x)\right)$ and
$h\left(f\left(\F_f(x)\right)\right)=g\left(h\left(\F_f(x)\right)\right)$
\end{theorem}

In \cite{hipush} they put the following:
\begin{problem}
If $f$ is partially hyperbolic and plaque expansive at $\F^c$, is
$\F^c$ the unique $f-$invariant foliation tangent to $E^c$?
\end{problem}

We heard the following problem from Charles Pugh:
\begin{problem}
Let $f:M\to M$ be a partially hyperbolic diffeomorphism and assume
that it has a central foliation by compact leaves. Is it true that
the volume of central leaves is uniformly bounded?, is it true that
the central foliation is plaque expansive?, is it true that there is
a fibration $p:M\to N$ whose fibers are the central leaves and an
Anosov diffeomorphism $g:N\to N$ such that $p$ is a semiconjugacy
between $f$ and $g$?.
\end{problem}
Maybe the last question goes a little beyond of what Charles Pugh
was asking, and maybe it goes a little beyond the reality.

\begin{subsubsection}{Lyapunov Stable}
Let $f:M\to M$ leave invariant a bundle $E\subset TM$. We say that
$f$ is Lyapunov stable in the direction $E$ if for any $\epsilon>0$
there is $\delta>0$ such that for any $C^1$ path $\gamma$ tangent to
$E$, $length(\gamma)<\delta$ implies
$length\left(f^n\gamma\right)<\epsilon$ for every $n\geq 0$.
\begin{theorem}\label{lyap}
Let $f:M\to M$ admit a splitting $TM=E^{cs}\oplus E^u$, where
$Df|E^u$ is expanding.
Let us assume that $f$ is Lyapunov stable in the direction $E^{cs}$.
Then $E^{cs}$ is tangent to a unique lamination $W^{cs}$. Moreover,
$W^{cs}$ is plaque expansive.
\end{theorem}

\begin{proof}
Let us assume for simplicity that $|Df^{-1}|E^u|\leq\mu^{-1}$,
$\mu>1$. The existence of the lamination tangent to $E^{cs}$ follows
directly from the technics in Theorem (7.5) of \cite{hipush}. Let us
see that $W^{cs}$ is plaque expansive. Let $\Pe$ be a plaquation of
$W^{cs}$ and let $\{x_n\}_n$ and $\{y_n\}_n$ be two $\nu$-pseudo
orbits respecting $\Pe$. Let us assume that $d(x_n,y_n)<\nu$ for
every $n\in\Z$. We want to show that if $\nu$ is small enough, then
$x_0$ and $y_0$ are in the same plaque. Let
$z_n=W_{loc}^{cs}(x_n)\cap W_{loc}^u(y_n)$. Then $z_n$ is an
$\delta$-pseudo orbit respecting $\Pe$ and $d(z_n,y_n)<\delta$ for
every $n\in\Z$ where $\delta$ goes to zero with $\nu$. If $z_0$ and
$y_0$ lie in a common plaque then $x_0$ and $y_0$ are in the same
plaque and we are done. Let us call $\sigma_k=\sup_{n\in\Z}
d\left(f^k(y_n),f^k(z_n)\right)$. There is $\epsilon_0$ independent
of $\delta$ and $\nu$ such that if $\sigma_j\leq\epsilon_0$ for
$0\leq j\leq k-1$ then $\sigma_k\geq\mu^k\sigma_0$. Take $k_0+1$ the
first time $\sigma_k>\epsilon_0$, which exists since we assume by
contradiction that $\sigma_0\neq 0$. Then
$\epsilon_0\geq\sigma_{k_0}\geq\frac{\epsilon_0}{|Df|E^u|}=\epsilon_1$.

Let us fix a small $\epsilon$ and let us take $\delta$ from the
definition of Lyapunov stable. Given $k\geq 0$, we have that
$\left\{f^k(y_n)\right\}_n$ and $\left\{f^k(z_n)\right\}_n$ are
$\epsilon$-pseudo orbits respecting $\Pe$. In fact
$d\left(f\left(f^k(y_n)\right),f^k(y_{n+1})\right)=d\left(f^k\left(f(y_n)\right),
f^k(y_{n+1})\right)<\epsilon$ since
$d\left(f(y_n),y_{n+1}\right)<\delta$. Let us call
$p_n=f^{k_0}(y_n)$ and $q_n=f^{k_0}(z_n)$. We have that
\begin{eqnarray*}
\mu d(p_n,q_n)&\leq& d\left(f(p_n),f(q_n)\right)\\
&\leq& d\left(f(p_n),p_{n+1}\right)+d\left(p_{n+1},q_{n+1}\right)+d\left(f(q_n),q_{n+1}\right)\\
&\leq&2\epsilon+\sigma_{k_0}
\end{eqnarray*}
Taking $n$ such that $d(p_n,q_n)\geq\sigma_{k_0}-\epsilon$, we get
that
$\epsilon_1\leq\sigma_{k_0}\leq\left(\frac{2+\mu}{\mu-1}\right)\epsilon$.
So, taking $\epsilon$ small, since $\epsilon_1$ is fixed we get a
contradiction and thus we get the theorem.
\end{proof}

\begin{corollary}
Let $f$ be a partially hyperbolic diffeomorphism. If $f$ and
$f^{-1}$ are Lyapunov stable in the direction $E^c$ then $E^c$ is
tangent to a unique lamination $W^c$. Moreover, $W^c$ is plaque
expansive. The same holds if $1/C\leq m(Df^n|E^c)\leq |Df^n|E^c|\leq
C$ for every $n\geq 0$ and some constant $C>0$ which is the case
when $Df|E^c$ is an isometry.
\end{corollary}
This last case appeared in \cite{hipush}, but no proof were
available since then.

\end{subsubsection}

\end{subsection}

\end{section}

\section{Robust transitivity}\label{rtransitivity}

A diffeomorphism $f$ of a closed manifold is $C^r$ robustly (stably)
transitive if it belongs to the $C^r$ interior of the transitive
diffeomorphisms.

Transitive Anosov diffeomorphisms are examples of such
diffeomorphisms (it is well known that the transitivity of any
Anosov diffeomorphism is an open question). The first nonhyperbolic
examples were given by Shub (\cite{sh2}) in $\mathbb{T}^4$. Later
Ma\~n\'e gave an example on $\mathbb{T}^3$ (\cite{mane1}), in
dimension $2$ robust transitivity implies Anosov \cite{mane2}. The
list of examples that are known to be robustly transitive is very
small:
\begin{itemize}
\item Transitive Anosov diffeomorphisms.
\item Some derived from Anosov (see \cite{mane1,sh2}).
\item Perturbations of $f\times Id$ where $f$ is a transitive Anosov
diffeomorphism $Id$ is the identity map of any closed manifold (see
\cite{bodi}).
\item Perturbations of the time-one map of a transitive Anosov flow
(see \cite{bodi}).
\item Examples that are not partially hyperbolic but
 presenting some form of weak hyperbolicity  (see \cite{bovi}).

\end{itemize}

An important tool in order to construct examples are the
center-stable blenders first introduced in \cite{bodi}. The blenders
resemble a  high dimensional skew horseshoe and their more important
property is that in their presence a quasi-transversal intersection
between stable and unstable manifolds of hyperbolic sets of
different indices turns out to be persistent under perturbations.
For instance, this enables Bonatti and D\'\i az to prove that the
closure of the stable manifold of certain periodic point contains a
stable manifold of greater dimension. This property, combined with
some global property of the original diffeomorphism (the presence of
two periodic points that persistently have dense stable or unstable
manifold), implies that the transverse homoclinic points of these
periodic points are dense giving the desired transitivity.

Although there exists this abstract construction of robustly
transitive diffeomorphisms, the first problem is to enlarge the
known set of examples. It seems that many of the examples that are
known to be stably ergodic should be robustly transitive. For
instance:

\begin{problem} Is the time-one map of the geodesic flow of
a surface of negative curvature robustly transitive? And the
partially hyperbolic automorphisms of three dimensional nilmanifolds
of section \ref{sacksteder}?
\end{problem}

Moreover, there are no satisfactory sufficient conditions implying
robust transitivity. There are some theorems with necessary
conditions in the $C^1$ category: some weak form of hyperbolicity is
needed (we shall explain it better below) and if  $f$ is partially
hyperbolic and the center bundle is one-dimensional,generically, at
least one of the strong foliations must be minimal (see
\cite{bodiur,rhrhur2}).

\begin{problem} Let $f$ be a partially hyperbolic diffeomorphism.
Does minimality of the strong unstable foliation imply that $f$ is
robustly transitive? And if, in addition, we demand $f$ to have the
accessibility property?
\end{problem}

In \cite{dipuur} and \cite{bodipu} it is proved that some amount of
hyperbolicity is needed in order to obtain robust transitivity (at
least in the $C^1$ topology). In fact, for surface diffeomorphisms
Ma\~n\'e results (see \cite{mane2}) implies that $C^1$ robustly
transitive diffeomorphisms are Anosov and, of course, there are not
robustly transitive diffeomorphisms of $\S^1$.

Let us explain the results of \cite{bodipu} that generalize to any
dimension those of \cite{dipuur} for dimension 3. A continuous
invariant splitting $TM=E\oplus F$ is called dominated if there
exists $n\in \mathbb{N}$ such that
$$||D_xf^n|_{E}||. ||D_{f^n(x)}f^{-n}|_{F}||<1/2.$$

Moreover, we will say that a continuous invariant splitting
$TM=E_1\oplus\dots \oplus E_k$ is dominated if, for each
$i=1,\dots,k$, the bundles $E=E_1\oplus\dots \oplus E_i$ and
$F=E_i\oplus\dots \oplus E_k$ define a dominated splitting in the
latter sense.

\begin{theorem}\cite{dipuur, bodipu} Let $f$ be a $C^1$ robustly transitive
diffeomorphism. Then, $f$ admits a dominated splitting
$TM=E_1\oplus\dots \oplus E_k$ such that the Jacobian of
$Df^n|_{E_1}$ and  $Df^{-n}|_{E_k}$ decrease exponentially with $n$.
\end{theorem}

This  theorem says that robustly transitive diffeomorphism are
partially ``volume"  hyperbolic.

\begin{problem} Is it possible to build a theory of
partially volume hyperbolic diffeomorphism? Observe that, in the
easier case of dimension 3, at least one of the bundles is
hyperbolic (it is one dimensional) but the other could be only
volume hyperbolic and, then, it is not known if it is integrable.
Are these strong bundles uniquely integrable as in the partially
hyperbolic case?
\end{problem}

This motivates also the following problem: Give sufficient
conditions (optimal) in such a way that a volume preserving
diffeomorphism on a $3$ dimensional manifold admitting a splitting
of the form $TM=E^{cs}\oplus E^u$ be ergodic. Here, $E^{cs}$ is
volume contracting and $E^u$ expands vectors. What about the
following?
\begin{problem}
Is it true that a $C^1$ open and $C^{\infty}$ dense set of volume
preserving diffeomorphism admitting a splitting $TM=E^{cs}\oplus
E^u$ is ergodic? Here again $E^{cs}$ is volume contracting and $E^u$
expands vectors. What about in dimension $3$?.
\end{problem}
Aside from the case some hypothesis is made on Lyapunov exponents,
very little is known to guaranty ergodicity. Maybe in dimension $3$,
if $E^{cs}$ cannot be split then its Lyapunov exponents are $C^r$
typically negative, $r\geq 2$.

Let us also mention that Horita and Tahzibi \cite{horta} and
independently Saghin have proven that stable ergodicity among $C^1$
symplectic diffeomorphisms implies partial hyperbolicity.

Related to transitivity, there is the following counterpart of
theorem \ref{accessimpliesweakerg}:
\begin{theorem}\cite{br1}
If $\Omega(f)=M$ and $GC(x)$ is dense for some point $x$ then $f$ is
transitive.
\end{theorem}
The reader may find also a proof of this in \cite{rhrhur2}. So we
have the following counterpart of Corollary
\ref{ergodictransversal}.
\begin{problem}
Let $f:M\to M$ be a partially hyperbolic diffeomorphism. If $P^c$ is
a compact $f-$invariant manifold tangent to the central direction,
$f|P^c$ is robustly transitive and $f$ has the stable accessibility
property, is $f$ robustly transitive?
\end{problem}
We expect that the answer to this problem is no. In fact if we take
the non transitive Anosov flow $\phi$ constructed in \cite{frwi}, it
has the stable accessibility property and we can take a time $t$ in
such a way that for some closed orbit $\phi_t|P^c$ is an irrational
rotation. Of course $\phi_t|P^c$ is not robustly transitive, but
maybe if we multiply by an Anosov diffeomorphism $A$ on $\T^2$ and
make a perturbation $g$ to get that $g|P^c\times\T^2$ is robustly
transitive... Compare with proposition (8.4) of \cite{hipush}.

\begin{section}{Classification}\label{clasif}
As we have already said, the problem of the classification of
partially hyperbolic systems is widely open. In \cite{push3} Charles
Pugh and Mike Shub posed the following problem related to the
topology of the manifold supporting a partially hyperbolic system:
\begin{problem}
If a manifold $M$ supports a partially hyperbolic diffeomorphism,
does $M$ fibers over a lower, positive dimensional manifold?
\end{problem}

\begin{subsection}{Transitive systems in dimension
3}\label{transbowi} In \cite{bowi} Christian Bonatti and Amie
Wilkinson have done some substantial  advances  in an attempt of
classification of  transitive partially hyperbolic diffeomorphisms
on three manifolds.  Their main hypothesis concern the existence of
an invariant embedded circle $\gamma$ (observe that $\gamma$ will
always be tangent to $E^c$) and the behavior of the invariant
foliations around it. In their own words, if a transitive partially
hyperbolic diffeomorphism $f$ looks like the perturbations of a skew
products or of the time-1 map of an Anosov flow in just a small
region of a 3-manifold, then $f$ is the perturbation of a skew
product or the time-1 map of an Anosov flow.

$W^s_\delta(\gamma)$ and $W^u_\delta(\gamma)$ will denote the union
of the strong stable and strong unstable segments, respectively, of
length $\delta$ through the points of $\gamma$.
 Let us state the theorems:

\begin{theorem} \label{skewbw} Let $f$ be a partially hyperbolic diffeomorphism of
a 3-manifold $M$. Assume that there is an embedded circle $ \gamma$
such that $f(\gamma)=\gamma$. Suppose there exists $\delta>0$ such
that $W^s_\delta (\gamma)\cap W^u_\delta(\gamma)\setminus \gamma$
contains a connected circle. Then:
\begin{enumerate}
\item  $f$ is dynamically coherent.
\item Each center leaf is a circle and the center foliation is a
Seifert bundle on M.
\item If the center-stable and the center-unstable foliations are
transversely orientable, then $M$ is a $\S^1$-bundle over $\T^2$,
and $f$ is conjugate to a (topological) skew product over a linear
Anosov map of $\T^2$.
\item If the center-stable or the center-unstable foliations are not
orientable, then a covering of $M$ corresponding to the possible
transverse orientations is a $\S^1$-bundle and the natural lift
$\tilde f$ of $f$ is conjugate to a (topological) skew product over
a linear Anosov map of $\T^2$.

\end{enumerate}
\end{theorem}

In this result they adopt a more general definition of a skew
product. They say that a homeomorphism $F$ of any circle bundle is a
skew product over an Anosov map $A$ if $F$ preserves the fibration
and projects to $A$. Observe that this definition includes the
examples of R. Sacksteder on nilmanifolds (see \cite{sa} and
subsection \ref{sacksteder}).

\begin{theorem}\label{anosvoflowbw}
Let $f$ be a partially hyperbolic dynamically coherent
diffeomorphism on a compact 3-manifold $M$. Assume that there is a
closed center leaf $\gamma$ which is periodic under $f$ and such
that each center leaf in $W^s_{loc}(\gamma)$ is periodic for f.

Then:
\begin{enumerate}
\item there is an $n\in \mathbb{N}$ such that $f^n$ sends every center
leaf to itself.
\item there is an $L>0$ such that for any $x\in M$ the length of the
smaller center segment joining $x$ to $f^n(x)$ is bounded by $L$.
\item each center-unstable leaf is a cylinder  or a plane (according
it contains a closed center leaf or not) and is trivially
bi-foliated by center and strong unstable leaves.
\item the center foliation supports a
transitive expansive continuous flow.
\end{enumerate}
\end{theorem}

\end{subsection}

\begin{subsection}{Growth of curves in dimension $3$.}\label{growthcurves}
The results in this subsection follows essentially from
\cite{brbuiv}. Their idea is to analyze the action of $f$ on
homology and thus get some restrictions on the homotopy type. Here
we push this argument a little more. The results here are for
dimension $3$. One of the main tools here is Novikov theorem on Reeb
components.

Given  a compact manifold $M$ and $x\in\tilde M$, let us define
$\vol_x(r)=vol\left(B(x,r)\right)$. Notice that there is $C>0$ such
that $\vol_x(r)\leq C\vol_y(r)$ for any two point $x$ and $y$. So
let us fix $x_0\in\tilde M$ and call $\vol(r)=\vol_{x_0}(r)$.
\begin{proposition}\label{volgro}
Let $f:M\to M$ be a partially hyperbolic diffeomorphism on a three
dimensional manifold. Assume that either $E^s\oplus E^u$ or
$E^c\oplus E^u$ is integrable. Then there is a a constant $C>0$ such
that if $I\subset\tilde M$ is an unstable arc then $length(I)\leq
C\vol\left(diam(I)\right)+C$. Moreover, if $x\in W^u(y)$ then
$d^u(x,y)\leq C\vol\left(d(x,y)\right)+C$.
\end{proposition}
Thus, proposition \ref{volgro} gives something that recalls the
quasi-isometric property used by Brin in \cite{br3}. In fact, in
some cases we can get something closer, for example for
nilmanifolds, as they have polynomial growth of volume we get:
\begin{corollary}\label{expgrowth}
In the setting of proposition \ref{volgro}, if $M$ is a nil-manifold
then there is $C>0$ such that $length(I)\leq C
\left(diam(I)\right)^4+C$, and in fact if $x\in W^u(y)$ then
$d^u(x,y)\leq Cd(x,y)^4+C$.
\end{corollary}
But if the manifold is the unit tangent bundle of a hyperbolic
surface, then the volume growth exponentially and hence proposition
\ref{volgro} only give that $d(f^n(x),f^n(y))$ growth linearly.

Another property that is useful is the following:
\begin{proposition}
Let $f:M\to M$ be a diffeomorphism on a three dimensional manifold,
and assume it leave invariant a codimension one torus $T$ such that
$f|T$ leave invariant an expanding foliation, then $T$ is not the
boundary of a solid torus.
\end{proposition}

This proposition together  the following, that is also a consequence
of Novikov's theorem, will give a good description of the homotopy
type of some partially hyperbolic systems on some three dimensional
manifolds.
\begin{proposition} Let $f:M\to M$ be a partially hyperbolic
diffeomorphism with $\dim M=3$. If either $E^s\oplus E^u$, $E^{cs}$
or $E^{cu}$ is integrable then $M$ is a $K(\pi_1(M),1)$ manifold and
thus its universal covering is contractible.
\end{proposition}

Using all this facts it is proven:
\begin{theorem}\cite{brbuiv}
Let $f:M\to M$ be a partially hyperbolic diffeomorphism with $\dim
M=3$. If either $E^s\oplus E^u$, $E^{cs}$ or $E^{cu}$ is integrable
and $\pi_1(M)$ is abelian then the action of $f$ on its first
homology group is partially hyperbolic.
\end{theorem}

The idea is the following, take a point $x_0\in\tilde{M}$, the
universal covering and let $\Gamma=\pi_1(M)$ acts on $\tilde{M}$ by
the deck transformations. Then $\gamma:\Gamma\to\tilde{M}$,
$\gamma(L)=L(x_0)$ is a quasi isometry, that is, there is $C>0$,
independent of $x_0$ such that $\frac{1}{C}d(\gamma(L_1),
\gamma(L_2))-C\leq d_{\Gamma}(L_1,L_2)\leq Cd(\gamma(L_1),
\gamma(L_2))+C$ for every $L_1$ and $L_2$, where
$d_{\Gamma}(L_1,L_2)$ is the word length. Then, if $\pi_1(M)$ is
nilpotent, it can be proven that there are $L_1, L_2\in\Gamma$ such
that $d_{\Gamma}(f_*^n(L_1),f_*^n(L_2))\geq\sigma^n$ for some
$\sigma>1$ for every $n\geq n_0$.

It would be interesting to go to higher dimensions, for example,
\begin{problem}
If $f$ is a partially hyperbolic diffeomorphism on $\T^n$, $n\geq
3$, is it homotopic to a partially hyperbolic automorphism?
\end{problem}
Also here is quite evident the importance of the integrability to
the classification problem. Let us finish with other type of
problem.
\begin{problem}
Let $E\subset TM$ be a plane field on a three dimensional manifold.
Assume that $E$ integrates to a lamination on a closed  set
$\Lambda\subset M$. Is it true Novikov's theorem in this setting?,
if $\eta$ is a homotopically trivial closed curve transversal to
$E$, and $\eta\cap\Lambda\neq\emptyset$, does $\Lambda$ contains a
torus tangent to $E$?
\end{problem}

In \cite{rhrhur3} we recall all this results to find the ergodic
partially hyperbolic diffeomorphisms on dimension $3$

\end{subsection}

\end{section}

\end{document}